\documentclass[10pt]{article}

\usepackage[english]{babel}
\usepackage{amsmath}
\usepackage{amsfonts}
\usepackage{amssymb}
\usepackage{amsthm}
\usepackage[all]{xy}

\newcommand{\CC}{\mathbb{C}}
\newcommand{\ZZ}{\mathbb{Z}}

\newcommand{\RR}{\mathbb{R}}

\newcommand{\XX}{\mathfrak{X}}

\newcommand{\Wedge}{\textstyle{\bigwedge}}

\newcommand{\la}{\left\langle}
\newcommand{\ra}{\right\rangle}
\newcommand{\To}{\longrightarrow}
\newcommand{\rra}{\rightrightarrows}

\renewcommand{\AA}{\mathcal{A}}
\newcommand{\FF}{\mathcal{F}}

\newcommand{\Section}[1]{\Section{#1} \setcounter{equation}{1}}

\newcommand{\fiber}[2]{\mathbin{{}_{#1}\times_{#2}}}
\newcommand{\inner}{\mathbin{\raise0.1ex\hbox{$\lrcorner$}}}

\DeclareMathOperator{\End}{End}
 
\DeclareMathOperator{\Hom}{Hom}

\DeclareMathOperator{\Aut}{Aut}

\DeclareMathOperator{\Pic}{Pic}
\DeclareMathOperator{\Ad}{Ad}
\DeclareMathOperator{\ad}{ad}
\DeclareMathOperator{\im}{im}

\DeclareMathOperator{\spin}{Spin}

\theoremstyle{plain}
\newtheorem{theorem}{Theorem}[section]
\newtheorem*{theorema}{Theorem}
\newtheorem{lemma}[theorem]{Lemma}
\newtheorem{proposition}[theorem]{Proposition}

\newtheorem{corollary}[theorem]{Corollary}

\theoremstyle{definition}
\newtheorem{remark}[theorem]{Remark}
\newtheorem{definition}[theorem]{Definition}
\newtheorem{example}[theorem]{Example}

\begin{document}

\title{Geometric quantization of Hamiltonian actions of Lie algebroids and Lie groupoids}
\date{}
\author{Rogier Bos\thanks{Institute for Mathematics, Astrophysics and Particle Physics; Radboud University Nijmegen, The Netherlands, \texttt{Email:\ R.Bos@math.ru.nl},
this work is part of the research programme of the `Stichting voor Fundamenteel Onderzoek der Materie (FOM)', 
which is financially supported by the `Nederlandse Organisatie voor Wetenschappelijk Onderzoek (NWO)}}

\maketitle

\begin{abstract}
We construct Hermitian representations of Lie algebroids and associated unitary representations of Lie groupoids by a geometric 
quantization procedure. For this purpose we introduce a new notion of Hamiltonian Lie algebroid actions.
The first step of our procedure consists of the construction of a prequantization line bundle. 
Next, we discuss a version of K\"{a}hler quantization suitable for this setting. 
We proceed by defining a Marsden-Weinstein quotient for our setting and prove a ``quantization commutes with reduction''
theorem. We explain how our geometric
quantization procedure relates to a possible orbit method for Lie groupoids. 
Our theory encompasses the geometric quantization of symplectic manifolds, 
Hamiltonian Lie algebra actions, actions of bundles of Lie groups,
and foliations, as well as some general constructions from differential geometry.
\\
\\
\textsc{\footnotesize{Keywords: Lie groupoid, Lie algebroid, 
prequantization, geometric quantization, orbit method.}}\\
\\
\small{\textit{2000 Mathematics Subject Classification}. Primary:22A22, 58H05, 53D50. Secondary: 57R30, 53D20.}
\end{abstract}
\newpage
\tableofcontents
\section*{Introduction}
The aim of this paper is to give a method to construct Hermitian representations of Lie algebroids and associated 
unitary representations of Lie groupoids.
An important way of constructing representations of Lie algebras and Lie groups is by geometric quantization (cf.\ e.g.\ \cite{Woo},
\cite{HH}, \cite{Kir}). In this paper this procedure will be generalized to Lie algebroids and Lie groupoids.

A groupoid is a small category in which all arrows are invertible. In particular, it consists of a base set $M$ and
a set of arrows $G$, with a number of structure maps: the source and target map $s,t:G\to M$, composition
$m:G\fiber{s}{t}G\to G$, an inverse map $i:G\to G$ and a unit map $u:M\to G$, 
satisfying certain properties. If there exist smooth  
structures on the sets $M$ and $G$ such that the structure maps behave well, then $G$ is
called a Lie groupoid. A general reference for Lie groupoids is the book by K. Mackenzie \cite{Mac}. Lie groupoids 
are useful models for orbifolds (cf.\ \cite{Pr}), orbit spaces of Lie group
actions, and foliations (cf.\ \cite{MM1}). They are also used in the study of manifolds with boundaries or corners 
(cf.\ \cite{MoNi}) and play an important r\^{o}le in Poisson geometry (cf.\ \cite{Mac}, \cite{Wei1}). 

In the previous examples the groupoid is seen as a model for a singular space. 
This paper takes a different perspective: groupoids model a generalized notion of symmetry.
One usually studies the symmetry of an object $X$, by studying its set of automorphisms $\Aut(X)$, which has the structure of a group.
The group $\Aut(X)$ is often very large and one instead studies morphisms $G\to\Aut(X)$ for smaller groups $G$, called actions of $G$.  
One can go one step further and study the symmetry of a map $f:X\to Y$. An automorphism of $f:X\to Y$ consists of automorphisms 
$\phi\in\Aut(X)$ and $\psi\in\Aut(Y)$ such that $\psi\circ f=f\circ\phi$. Note that if $f$ is surjective, the automorphism $\phi$ of $X$ fixes the automorphism $\psi$ of $Y$, hence the automorphisms of $f$ form a subgroup of the automorphisms of $X$. If $Y$ is a set, then such an automorphism of a map actually consists of a family of isomorphisms
$\{\phi_{y}:f^{-1}(y)\to f^{-1}(\psi(y))\}$. If $\psi(y)=y$, then $\phi_{y}$ is called an internal symmetry of the map, else it is called an external symmetry (cf.\ \cite{Wei2}).
The union of all internal and external symmetries has the structure of a groupoid $\Aut(f)\rra Y$. For example, the symmetry of a principal 
$H$-bundle $f:P\to M$ (for a Lie group $H$ and a manifold $M$) is described by the gauge groupoid $P\times_{H}P\rra M$.
In this paper we shall study morphisms of groupoids $G\to\Aut(f)$, called actions of groupoid. In particular, we shall construct linear actions (representations) from Hamiltonian actions.

The infinitesimal structure associated to a Lie groupoid is a Lie algebroid. A Lie algebroid is a smooth vector bundle
$\AA\to M$ with a Lie bracket on the space of smooth sections $\Gamma^{\infty}(\AA)$ and a bundle morphism $\rho:\AA\to T M$, called the
anchor, satisfying a Leibniz identity for $f\in C^{\infty}(M)$ and $\tau,\sigma\in\Gamma^{\infty}(\AA)$, viz.\
\[[\tau,f\sigma]=f[\tau,\sigma]+(\rho(\tau)\cdot f)\,\sigma.\]

We shall often assume that Lie algebroids $\AA$ are regular, i.e.\ $\im(\rho))\subset T M$ has 
locally constant rank. If $\AA$ is the Lie algebroid associated to a Lie groupoid $G$, then this condition 
implies that the orbit foliation on $M$ of the groupoid $G$ is regular.
A Lie groupoid with regular orbit foliation is called regular. For example, transitive Lie groupoids, \'{e}tale Lie groupoids
and smooth bundles of Lie groups are regular. The regularity assumption is necessary to give proofs of some of the statements, 
but many constructions are possible to some extent in singular cases too. 

We now give an outline of the paper, including some more details on the content. In the first section we recall 
the notion of a Lie groupoid or Lie algebroid action on a map $J:S\to M$. 
The introduction of Hamiltonian actions of Lie algebroids (and Lie groupoids) proceeds in two steps.
First, for the case that $J:S\to M$ is a surjective submersion endowed with a family of 
symplectic forms $\omega$ we introduce the notion of an internally Hamiltonian action of a Lie algebroid. 
The word ``internal'' refers to the fact that we only consider the symmetry of the fibers $J^{-1}(m)$ for $m\in M$, 
which is represented by the action of the isotropy Lie algebra's $\AA_{m}$ for $m\in M$.
This action is internally Hamiltonian if there exists an internal momentum
map, which is a map
$\mu:S\to J^{*}\ker(\rho)^{*}$, satisfying certain natural conditions.

The second step considers the case that $\omega$ extends to a closed form $\tilde{\omega}$ on $S$ (which we shall call a $J$-presymplectic form).
Then, one can proceed by defining the notion of a Hamiltonian action, as is done Section 1.4. An action will be
called Hamiltonian if there exists a momentum map 
\[\tilde{\mu}:S\to J^{*}\AA^{*},\]
satisfying natural conditions. 
We shall give many examples to motivate this definition. Some of the examples will return throughout the paper.

Section two is devoted to the construction of prequantization line bundles with a representation of the Lie
algebroid, based on the data of a Hamiltonian Lie algebroid action. We introduce longitudinal \v{C}ech cohomology to study
such line bundles endowed with a connection. The main result of Section 2 is, 
summarizing Theorems \ref{exline}
and \ref{preqrep},
\begin{theorema}
If a Lie algebroid $\AA$ acts in a Hamiltonian fashion on $(J:S\to M,\tilde{\omega})$ and \mbox{$[\tilde{\omega}]\in H^{J,dR}(S)$} 
is integral, then there exists a prequantization line bundle carrying a Hermitian representation of $\AA$.
\end{theorema}

In Section 2.4 we briefly
discuss the possible integrability of such a representation to a representation of an integrating Lie groupoid for the Lie algebroid.

In the third section we obtain a representation of the Lie algebroid through generalized K\"{a}hler
quantization. To this effect we need $J:S\to M$ to be a bundle of compact K\"{a}hler manifolds. 
The main result is (cf. Theorem \ref{qu})
\begin{theorema}
If a Lie algebroid $\AA$ acts in a Hamiltonian fashion on $(J:S\to M,\tilde{\omega})$,\\  
\mbox{$[\tilde{\omega}]\in H^{J,dR}(S)$} is integral and $J:S\to M$ is a bundle of K\"{a}hler 
manifolds, then there exists a geometric quantization $Q(S,\tilde{\omega})$ carrying a Hermitian 
representation of $\AA$.
\end{theorema}
Next, we study the symplectic reduction of Hamiltonian groupoid actions (a generalized Marsden-Weinstein quotient).
We introduce an internal quotient $(I_{G}\backslash\mu^{-1}(0_{M}),\omega_{0})$ and a `full' quotient 
$(G\backslash\mu^{-1}(0_{M}),\tilde{\omega}_{0})$.
We also introduce internal quantum reduction $Q(S,\tilde{\omega})^{I_{G}}$ and full quantum reduction
$Q(S,\tilde{\omega})^{G}$.

Finally we prove a ``quantization commutes with reduction theorem'' for regular proper Lie groupoids,
\begin{theorema} (cf. Theorem \ref{qure} and Corollary \ref{qure2})
If $G$ is a proper groupoid acting in a proper, free and Hamiltonian fashion on a bundle of K\"{a}hler manifolds $(J:S\to M,\tilde{\omega})$, 
and $[\tilde{\omega}]\in H^{J,dR}(S)$ is integral, then there exist
isomorphisms of vector bundles
\[Q(I_{G}\backslash\mu^{-1}(0_{M}),\omega^{0})\overset{\simeq}{\to}(Q(S,\tilde{\omega}))^{I_{G}}\]
and
\[Q(G\backslash\mu^{-1}(0_{M}),\tilde{\omega}^{0})\overset{\simeq}{\to}(Q(S,\tilde{\omega}))^{G}.\]
\end{theorema}
The proof strongly relies on the `quantization commutes with reduction theorem'' for compact Lie groups.

The orbit method as developed by Kirillov (cf.\ \cite{Kir}) is based on the idea that there should be a
certain correspondence 
between the irreducible unitary representations of a Lie group and the coadjoint orbits in the dual of its Lie algebra.
This method works very well for nilpotent Lie groups (cf.\ \cite{CG}) and  compact Lie groups (the Borel-Weil theorem). There are also nice 
results for reductive Lie groups (cf.\ \cite{Vog})
and even for quantum groups (cf.\ \cite{Kir2}). One might wonder if such a principle is also useful for Lie groupoids. In this paper we shall see
that the answer is affirmative, although a smooth family of coadjoint orbits is not the only ingredient to construct a 
representation. One needs some more structure to take care of the global topology. Moreover, one should realize that 
the coadjoint orbits are submanifolds of the dual of the Lie algebroid of the {\em isotropy} groupoid (which equals the dual
of the kernel of the anchor). Although the isotropy groupoid is in general not smooth, it plays an essential
r\^{o}le in understanding the representation theory of $G$.

The theory presented in this paper should be distinguished from the theory of symplectic groupoids and
their prequantization (cf.\ \cite{WX}). Symplectic groupoids were introduced by
Alan Weinstein and others in a program to geometrically quantize Poisson manifolds. This is not the purpose of this paper.
We neither assume any (quasi-)(pre-)symplectic structure on the Lie groupoid, nor do we construct the
geometric quantization of a Poisson manifold. Also, our notion of momentum map differs from the notion
in \cite{MW}.

The author would like to thank Eli Hawkins, Peter Hochs, Klaas Landsman, Alan Weinstein,  
Marius Crainic, Gert Heckman, Ieke Moerdijk, Michael Mueger and Hessel Posthuma 
for discussions, suggestions, interest and/or support at various stages of the research.

\section{Hamiltonian Lie algebroid actions}
\subsection{Actions of groupoids and Lie algebroids.}\label{act}
The material in this section is standard (see \cite{Mac}), except for the introduction of internally symplectic and
$J$-presymplectic actions of Lie groupoids and Lie algebroids.

Suppose $G\rra M$ is a Lie groupoid, with source map $s:G\to M$, target $t:G\to M$, unit map $u:M\to G$, 
inversion $i:G\to G$ and composition (or multiplication) \mbox{$m:G^{(2)}:=G\fiber{s}{t}G\to G$}. We shall use the notation
$i(g)=g^{-1}$, $m(g,g')=g\,g'$ and $u(m)=1_{m}$. We shall assume throughout that $G\rra M$ is source-connected, which means that all the fibers of $s$ are connected. 
Suppose $N$ is a smooth manifold and $J:N\to M$ a smooth map. 
\begin{definition}
A \textbf{smooth left action of $G$ on} $J:N\to M$ is a smooth map
\[\alpha: G\fiber{s}{J}N\to N\]
satisfying
\begin{align}
J(g\cdot n)&=t(g)\mbox{ for all }(g,n)\in G\fiber{s}{J}N,\\
1_{J(n)}\cdot n&=n\mbox{ for all }n\in N,\\
g\cdot(g'\cdot n)&=(g g')\cdot n\mbox{ for all }(g,g')\in G^{(2)}\mbox{ and }n\in J^{-1}(s(g')),
\end{align}
using the notation $g\cdot n:=\alpha(g,n)$.
\end{definition}
We shall also use the notation $\alpha(g)(n):=\alpha(g,n)$
There is an analogous notion of a right smooth action.
One can show that the $G$-orbits in $N$ are smooth submanifolds of $N$. Note that these orbits are equal to the fibers of the map 
\[\tilde{J}:=p\circ J:N\to M/G,\] 
where $p:M\to M/G$ is the quotient map from $M$ to the orbit space $M/G$. These orbits form
a regular foliation of $N$, if $G$ is a regular Lie groupoid.

\begin{example}
We mention three basic examples here, all arising from a groupoid $G\rra M$ itself. Firstly, one has the action 
of $G$ on the identity map $M\to M$ by $g\cdot s(g)=t(g)$. Secondly, one has the action of $G$ on $t:G\to M$ by multiplication 
$g\cdot g':=g\, g'$. Thirdly, one has a (in general non-smooth) action of $G$ on the associated isotropy groupoid (which is
in general not a smooth manifold, unless $G$ is regular)
\[I_{G}:=\{s^{-1}(m)\cap t^{-1}(m)\}_{m\in M}\to M\]
by conjugation $c(g)g':=g\,g'\,g^{-1}$.
\end{example}

Associated to a smooth groupoid action of $G$ on $J:N\to M$ is an \textbf{action Lie groupoid} $G\ltimes J$ over $N$. 
Its space of arrows 
is given by $G\fiber{s}{J}N$, the source map by $s(g,n):=n$, target map by $t(g,n)=g\cdot n$, multiplication by
$(h,g\cdot n)(g,n):=(h g,n)$ and inversion by $i(g,n):=(g^{-1},g\cdot n)$.

Suppose $J:S\to M$ is a smooth surjective submersion.
The vector bundle $\ker(T J)\subset T S$ is the integrable distribution underlying the foliation $\FF:=\{J^{-1}(m)\}_{m\in M}$ of $S$.
We shall use the notation $T^{J}S:=\ker(T J)$, $T^{*,J}S=:\ker(T J)^{*}$, $\XX_{J}^{\infty}(S):=\Gamma^{\infty}(\ker(T J))$
and $\Omega^{n}_{J}(S):=\Gamma^{\infty}(\Wedge^{n}\ker(T J)^{*}))$.
Moreover, there is an obvious differential $d^{J}:\Omega^{n}_{J}(S)\to \Omega^{n+1}_{J}(S)$, which gives
rise to a generalized de Rham cohomology denoted by $H^{n}_{J,dR}(S)$.

Suppose $\alpha$ is an action of $G$ on a smooth family of symplectic manifolds \mbox{$(J:S\to M,\omega)$}, where 
$\omega\in\Omega^{2}_{J}(S)$ is a smooth family of symplectic forms. 
The action is \textbf{internally symplectic} if it preserves the symplectic forms in the sense that
\[\alpha(g)^{*}\omega_{\sigma}=\omega_{g^{-1}\cdot\sigma},\]
for all $\sigma\in S$ and $g\in G_{J(\sigma)}^{J(\sigma)}$. 
This is just a ``family version'' of symplectic actions in the usual sense.

Suppose $\omega\in\Omega_{J}^{2}(S)$ extends to a closed 2-form $\tilde{\omega}\in\Omega^{2}_{\tilde{J}}(S)$. We call
a closed form \mbox{$\tilde{\omega}\in\Omega_{\tilde{J}}^{2}(S)$} that restricts to a smooth family of symplectic forms $\omega\in\Omega_{J}^{2}(S)$ a \textbf{$J$-presymplectic form}.
Note that
\[(\alpha(\gamma)^{*}\tilde{\omega}_{\gamma(m)\cdot \sigma})|_{T^{J}S}=\alpha(\gamma)^{*}\omega_{\gamma(m)\cdot \sigma},\]
for all open sets $U\subset M$, local bisections $\gamma:U \to G$, $m\in U$ and $\sigma\in S_{m}:=J^{-1}(m)$, since the local diffeomorphism
$\alpha(\gamma)$ maps $J$-fibers to $J$-fibers.
A local bisection is a map $\gamma:U \to G$ such that $s\circ\gamma=id|_{U}$ and $t\circ\gamma$ is a diffeomorphism onto its image.
The action is said to be \textbf{$J$-presymplectic} if  
\[\alpha(\gamma)^{*}\omega_{\gamma(m)\cdot \sigma}-\omega_{\sigma}=0,\]
for all open sets $U\subset M$, local bisections $\gamma:U \to G$, $m\in U$ and $\sigma\in S_{m}:=J^{-1}(m)$.
This is equivalent to 
\[\alpha(g)^{*}\omega_{g\cdot \sigma}=\omega_{\sigma},\]
for all $g\in G$ and $\sigma\in J^{-1}(s(g))$.

\begin{definition}\label{algac}
An \textbf{action of a Lie algebroid $(\pi:\AA\to M,\rho)$ on a map} $J:N\to M$ is a map
\[\alpha:\Gamma^{\infty}(\AA)\to\XX^{\infty}(N)\]
satisfying
\begin{align}
\alpha(X+Y)&=\alpha(X)+\alpha(Y)\\
\alpha(f X)&=J^{*}f\alpha(X)\\
\lbrack\alpha(X),\alpha(Y)]&=\alpha([X,Y])\label{eqnalg}\\
T J(\alpha(X))&=\rho(X)
\end{align}
for all $X,Y\in\Gamma^{\infty}(\AA)$ and $f\in C^{\infty}(M)$.
\end{definition}
\begin{example}
Any action of a Lie algebra $\mathfrak{g}$ on a manifold $N$ is a Lie algebroid action of $\mathfrak{g}\to *$ on 
$N\to *$.
\end{example}
\begin{example}
Every Lie algebroid $(\AA\to M,[\cdot,\cdot],\rho)$ acts on \mbox{$J:=id:M\to M$} via the anchor $\rho:\Gamma^{\infty}(\AA)\to\XX^{\infty}(M)$. 
\end{example}
\begin{example}
A Lie groupoid action $\alpha$ on a smooth map $J:N\to M$ gives rise to an action $\alpha'$ of the Lie algebroid $\AA(G)$ on 
$J:N\to M$ by
\[\alpha'(X)(n):=\frac{d}{d\tau}|_{\tau=0}\exp(\tau X)_{J(n)}n.\]
\end{example}

Suppose $(\pi:\AA\to M,\rho)$ is a regular Lie algebroid. 
An action of $\AA$ on a smooth family of symplectic manifolds $(J:S\to M,\omega)$ is 
\textbf{internally symplectic} if $\mathcal{L}_{\alpha(X)}\omega=0$ for all $X\in\Gamma^{\infty}(\ker(\rho))$.
We shall see some examples of internally symplectic actions in Section \ref{intham}.

The image $\rho(\AA)\subset T M$ of the anchor is an integrable distribution, which induces a foliation 
$\mathcal{F}_{\rho}$ on $M$.
Suppose $\AA$ acts on a smooth fiber bundle $J:S\to M$. Denote the projection of $S$ on the leaf space $M/\mathcal{F}_{\rho}$ by 
$\tilde{J}:S\to M/\mathcal{F_{\rho}}$.
Suppose $\tilde{\omega}\in\Omega_{\tilde{J}}^{2}(S)$ is a $J$-presymplectic form.
Note that $(\mathcal{L}_{\alpha(X)}\tilde{\omega})|_{T^{J}S}=\mathcal{L}_{\alpha(X)}\omega$ for all $X\in\Gamma^{\infty}(\AA)$
Then the action is said to be \textbf{$J$-presymplectic} if  $\mathcal{L}_{\alpha(X)}\omega=0$ for all 
$X\in\Gamma^{\infty}(\AA)$. Note that the action being $J$-presymplectic implies it being internally symplectic. We shall see many examples of $J$-presymplectic actions in Section
\ref{Hamacmommap}, since a Hamiltonian action as defined in this section is automatically $J$-presymplectic.

Associated to a Lie algebroid action of $(\AA\to M, \rho,[\cdot,\cdot])$ on $J:N\to M$ there is an \textbf{action Lie algebroid} $\AA\ltimes J$.
Denote the pullback of $\AA\to M$ along $J:N\to M$ by $J^{*}\AA\to N$. The space of sections $\Gamma^{\infty}(J^{*}\AA)$ is
generated as a $C^{\infty}(N)$-module by sections of the form $\,J^{*}X$ for $X\in\Gamma^{\infty}(\AA)$. 
A Lie bracket on the smooth sections is defined by
\begin{align*}[f\,J^{*}X,g\,J^{*}Y]&:=f\, g J^{*}[X,Y]+f\, (\alpha(X)\cdot g) J^{*}Y-g\,(\alpha(Y)\cdot f)J^{*} X,
\end{align*}
where $f,g\in C^{\infty}(N)$ and $X,Y\in\Gamma^{\infty}(\AA)$
and the anchor 
\[\rho':\Gamma^{\infty}(\AA\ltimes J)\to \XX^{\infty}(N)\]
is given by
\[\rho'(f J^{*} X):=f\alpha(X).\]

Suppose a Lie groupoid acts on a map $J:N\to M$. It induces an action of the Lie algebroid $\AA(G)$ on $J:N\to M$ and
the action Lie algebroid $\AA(G)\ltimes J$ is isomorphic to the Lie algebroid $\AA(G\ltimes J)$ associated to the action Lie groupoid. 

\subsection{Internally Hamiltonian actions and internal momentum maps.}
In this section we introduce the notion of internally weakly Hamiltonian Lie algebroid action.
This notion and the notion of internally strongly Hamiltonian Lie algebroid action, introduced in the next section,
should be seen as an intermediate stage towards defining Hamiltonian actions. They are separately treated 
for clarity and for their r\^{o}le in the orbit method. Examples of internal Hamiltonian actions are postponed to the next section.

Note that the isotropy groupoid $I_{G}$ of a (non-regular) Lie groupoid is not a smooth manifold. But for any $G$-orbit 
$G m\subset M$ the restriction $I_{G}|_{G m}$ is a smooth manifold. Hence $I_{G}\to M$ is a continuous family of smooth manifolds
in the subspace topology, i.e.\ a surjective continuous map of topological spaces such that each fiber is a smooth manifold with the subspace topology. 
\begin{example}
Consider the action of the circle $\mathbb{S}^{1}$ on the real plane $\RR^{2}$ by rotation. Consider the action groupoid 
$G=\mathbb{S}^{1}\ltimes\RR^{2}\rra \RR^{2}$. The isotropy groupoid is a continuous family of Lie groups with fiber $\mathbb{S}^{1}$ at 
$(0,0)$ and zero fiber elsewhere.
\end{example}
If $G$ is a regular Lie groupoid, then $I_{G}$ is a smooth family of Lie groups, i.e. $I_{G}\to M$ is a smooth family of manifolds
and each fiber has a Lie group structure smoothly depending on $m\in M$.
Let $\pi:\AA(I_{G})\to M$ be the smooth family of Lie algebras associated to $I_{G}$. It is naturally isomorphic to the 
kernel $\ker(\rho)$ of the anchor $\rho:\AA(G)\to T M$ of the Lie algebroid of $G$.

Suppose that  $(\pi:\AA\to M,\rho)$ is a regular Lie algebroid that acts on a smooth family of symplectic manifolds $(J:S\to M,\omega)$.
Denote the action by $\alpha:\Gamma^{\infty}(\AA)\to\XX_{\tilde{J}}^{\infty}(S)$.
Suppose that the action of $\AA$ is internally symplectic. Then $\alpha(X)\inner\,\omega$ is closed, i.e. 
\[d^{J}(\alpha(X)\inner\,\omega)=0\]
for all $X\in\Gamma^{\infty}(\ker(\rho))$.
Indeed, this follows from the Cartan homotopy formula
\[\mathcal{L_{\alpha(X)}}\omega=d^{J}(\alpha(X)\inner\,\omega )+\alpha(X)\inner\,d^{J}\omega,\]
in which the last term is zero, since $\omega$ is symplectic on $S$.
\begin{definition}
An internally symplectic action of a regular Lie algebroid $(\pi:\AA\to M,\rho)$ on a smooth family of symplectic manifolds 
$(J:S\to M,\omega)$ is called \textbf{internally weakly Hamiltonian} if there exists a smooth map 
$\mu:S\to\ker(\rho)^{*}$, such that
\[\xymatrix{S\ar[r]^-{\mu}\ar[d]^-{J}&\ker(\rho)^{*}\ar[dl]^{p}\\M&}\]
commutes and
\[d^{J}\la\mu,J^{*}X\ra=-\alpha(X)\inner\,\omega,\]
for all $X\in\Gamma^{\infty}(\ker(\rho))$. The map $\mu$ is called an \textbf{internal momentum map} for the $\AA$-action.
\end{definition}
\begin{remark}
One might view $\mu$ as a section in $\Gamma^{\infty}(J^{*}\ker(\rho)^{*})$.
\end{remark}
\begin{remark}
One should think of $[\alpha(X)\inner\,\omega]$ as a cohomological obstruction to the existence of a momentum map.
One has the following diagram
\[\xymatrix{C^{\infty}_{J}(S)\ar[r]^-{d^{J}_{0}}&\Gamma^{\infty}(\ker(d^{J}_{1}))\ar[r]&H^{1}_{J, dR}(S)\\&
\Gamma^{\infty}(\ker(\rho))\ar[u]^{-\alpha\inner\omega}\ar@{-->}[ul]^{\mu}\ar[ur]&}\]
where $\alpha\inner\omega(X):=\alpha(X)\inner\omega$ and the right diagonal arrow denotes  
the induced map on the quotient space. The vanishing of this map is a necessary condition for $\alpha\inner\omega$ to lift to a
map $\mu$. 
\end{remark}
\begin{definition}
A symplectic Lie groupoid action is \textbf{internally weakly Hamiltonian} if the associated Lie algebroid action is 
internally weakly Hamiltonian.
\end{definition}
Before we give examples of such actions we shall introduce the notion of internally strongly 
Hamiltonian actions in the next section.

\subsection{The coadjoint action and internal momentum maps.}\label{intham}
In this section we introduce the notion of internally strongly Hamiltonian actions and treat several examples.

Let $G$ be a regular Lie groupoid and $I_{G}$ the associated isotropy Lie groupoid. 
Recall that $G$ acts (from the left) on $I_{G}\to M$ by conjugation $G\fiber{s}{p}I_{G}\to I_{G}$,
\[c(g)g':=g g' g^{-1}.\]

The action by conjugation induces an action of $G$ on the smooth family of Lie algebras 
$\ker(\rho)\simeq\AA(I_{G})$ by
\[\Ad(g)X=\left.\frac{d}{d \tau}\right|_{\tau=0}c(g)\exp(\tau X),\]
where $X\in\AA(I_{G})_{m}$ and $g\in G_{m}$ for any $m\in M$.
This action is called the \textbf{adjoint action of $G$} and is the generalization of the adjoint action for Lie groups.

In turn, this induces the \textbf{adjoint action} of the Lie algebroid $\AA(G)$ on $\AA(I_{G})\to M$.
\[\ad(X)Y=\left.\frac{d}{d \tau}\right|_{\tau=0}\Ad(\exp(\tau X))Y,\]
where $X\in\AA(G)_{m}$ and $Y\in\AA(I_{G})_{m}$.
Note that for $X\in\Gamma^{\infty}(\AA(G))$ and $Y\in\Gamma^{\infty}(\AA(I_{G}))$
\[\ad(X)Y=[X,Y].\]

\begin{example}
A simple example is the pair groupoid $G=M\times M\rra M$. Conjugation is given by $c((m,n),(n,n))=(m,m)$. The kernel of the 
anchor is the zero bundle hence $\Ad$ is trivial on the fibers and $\ad: T M\to M\times\{0\}$ is the zero map.
\end{example}
\begin{example}
If $G$ is a Lie group, then $\Ad$ and $\ad$ coincide with the usual notions.
\end{example}

\begin{remark}
There is also the notion of an action up to homotopy (cf.\ \cite{ELW}). It turns out that the map $\ad(X)Y:=[X,Y]$ defines
an action up to homotopy of $\AA$ on itself. We shall not use this structure in our paper.
\end{remark}

One defines the \textbf{coadjoint action of $G$} on the dual bundle $\AA^{*}(I_{G})$ by
\[\la\Ad^{*}(g)\xi,X\ra:=\la\xi,\Ad(g^{-1})X\ra,\]
where $\xi\in\AA^{*}(I_{G})_{m}$ and $g\in G_{m}$.
Analogously, one defines the \textbf{coadjoint action of $\AA(G)$} on $\AA^{*}(I_{G})$ by
\[\la\ad^{*}(X)\xi,Y\ra:=\la\xi,\ad(-X)Y\ra,\]
which is obtained as the tangent map of $\Ad^{*}$.

\begin{definition}\label{locstr}
An internally weakly Hamiltonian action of a regular Lie algebroid $\AA$ on a smooth family of symplectic manifolds
$(J:S\to M,\omega)$ is \textbf{internally (strongly) Hamiltonian} if the momentum map $\mu:S\to\ker(\rho)^{*}$ is $\AA$-equivariant with respect to the 
coadjoint action of $\AA$ on $\ker(\rho)^{*}$, i.e.
\[\alpha(X)\cdot\la\mu,Y\ra=\la\ad^{*}(X)\,\mu,Y\ra.\]
\end{definition}
\begin{definition}
An internally weakly Hamiltonian Lie groupoid action is \textbf{internally (strongly) Hamiltonian} if the momentum map 
$\mu:S\to\ker(\rho)^{*}$ is $G$-equivariant with respect to the coadjoint action of $G$ on $\ker(\rho)^{*}$, i.e.
\[\mu(g\cdot\sigma)=\Ad^{*}(g)\cdot\mu(\sigma).\]
\end{definition}
\begin{example}\label{bundlelie}
One can consider for example smooth families of Lie algebras $\pi:\mathfrak{g}\to M$; 
in particular, a bundle of Lie algebras $P\times_{H}\mathfrak{h}$, where $H$ is a Lie group and $P\to M$ a principal $H$-bundle and the action of $H$ on $\mathfrak{h}$ is the adjoint action. More about internally Hamiltonian actions of such bundles can be derived from Example \ref{princ}.

In general one can remark the following. Suppose $\mathfrak{g}$ acts on a smooth family of symplectic manifolds $S:=\bigcup_{m\in M}S_{m}\to M$. 
Then a momentum map is a smooth map $S\to\mathfrak{g}^{*}$ that restricts to a momentum map in the classical 
sense on each fiber. For example, a smooth family of coadjoint orbits $S:=\{\mathcal{O}_{m}\subset\mathfrak{g}^{*}_{m}\}_{m\in M}$ carries a Hamiltonian action (namely the coadjoint action). The inclusion $S\hookrightarrow\mathfrak{g}^{*}$ is an internal momentum map.
\end{example}

\begin{example}\label{coadj}
Suppose $G\rra M$ is a regular Lie groupoid with associated Lie algebroid $(\AA\to M,\rho)$. Consider a family of coadjoint orbits 
\[\{\mathcal{O}_{m G}\}_{m G\in M/G}\]
in the dual of the kernel of the anchor $\ker(\rho)^{*}$. Suppose they form a smooth family
\[S:=\bigcup_{m G\in M/G}\mathcal{O}_{m G}\to M.\]
Then it has a symplectic structure at $m\in M$ given by the standard symplectic form 
on a coadjoint orbit $\mathcal{O}_{m G}\cap \ker(\rho)^{*}_{m}$ in the dual of the Lie algebra $\ker(\rho)_{m}$. The inclusion 
\[S:=\bigcup_{m G\in M/G}\mathcal{O}_{m G}\hookrightarrow\ker(\rho)^{*}\] 
is an internal momentum map for the coadjoint action on $S$ which is therefore internally Hamiltonian. 
This is an important observation concerning the orbit method for Lie groupoids. 
We shall come back to it later.
\end{example}

\begin{example}\label{princ}
Suppose $H$ is a Lie group and $\pi:P\to M$ a principal $H$-bundle. Denote the action of $H$ on $P$ by $\alpha$.
Suppose $H$ acts on a symplectic manifold $(S,\omega^{S})$ in a Hamiltonian fashion with momentum map $\mu:S\to\mathfrak{h}^{*}$.
Denote the action of $H$ on $S$ by $\beta$. 

Let $G\rra M$ denote the gauge groupoid $P\times_{H}P\rra M$.
Define a smooth bundle of smooth manifolds by
\[S':=P\times_{H}S.\]
The map $\bar{\pi}:[p,\sigma]\mapsto\pi(p)$ is well defined from $S'$ to $M$ and gives the bundle structure.
The following observations and lemma will be necessary to endow $S$ with the structure of a bundle of symplectic 
manifolds.

Note that, since $P$ is a principal $H$-bundle, the infinitesimal action
\[\alpha:P\times\mathfrak{h}\to T^{\pi}P\]
is an isomorphism of smooth vector bundles.
Moreover, it is equivariant with respect to the adjoint action of $H$ on $\mathfrak{h}$, hence it induces a diffeomorphism
\[\bar{\alpha}:P\times_{H}\mathfrak{h}\to T^{\pi}P/H.\]

\begin{lemma}\label{quotlem}
Suppose a Lie group $H$ acts properly and freely on a manifold $N$. Denote the action by $\gamma$.
Then 
\[T(N/H)\simeq (T N)/\mathord\sim,\] 
where the equivalence relation is generated by
\[h v\sim v\]
for all $v\in T N$ and $h\in H$ and
\[\gamma(X)_{n}\sim 0,\]
for all $X\in\mathfrak{h}$, $n\in N$, and where $\gamma:\mathfrak{h}\to\XX^{\infty}(N)$ denotes the infinitesimal action.
\end{lemma}
\begin{proof}
Consider the tangent map
\[T N\to T(N/H)\]
of the quotient map $N\to N/H$. It is surjective and the kernel is spanned by the elements mentioned above as 
one easily checks.
\end{proof}
\begin{remark}
One should compare this lemma to the fact that $H$ acts on $T^{*}N$ in Hamiltonian fashion 
with ``classical'' momentum map $\mu:T^{*}N\to\mathfrak{h}^{*}$, and that the Marsden-Weinstein 
quotient satisfies
\[T^{*}(N/H)\simeq\mu^{-1}(0)/H\]
(cf.\ e.g., \cite{HH}).
\end{remark}
Applying this lemma to $T S'=T (P\times_{H}S)$ one obtains
\[T(P\times_{H}S)\simeq (T P\times T S)/\mathord\sim,\]
and restricting to the vertical tangent space one has
\[T^{\bar{\pi}}(P\times_{H}S)\simeq (T^{\pi}P\times T S)/\mathord\sim.\]

The map $\alpha$ induces an isomorphism
\[(T^{\pi}P\times T S)/\mathord\sim\To(P\times\mathfrak{h}\times T S)/\mathord\sim_{1}\]
with the equivalence relation generated by 

\begin{align*}[p,X,\beta'(X)]&\sim_{1}[p,0,0]\\
(h\cdot p,h\cdot X, h\cdot v)&\sim_{1}(p,X,v),
\end{align*}
for all $p\in P$, $h\in H$, $X\in\mathfrak{h}$.
The map
\[(P\times\mathfrak{h}\times T S)/\mathord\sim_{1}\To P\times_{H}T S\]
given by $[p,X,v]\to [p,v-\beta'(X)]$ is again an isomorphism of smooth vector bundles. So we conclude that
\[T^{\bar{\pi}}(P\times_{H}S)\simeq P\times_{H}T S.\]

We can define a structure of a smooth bundle of symplectic manifolds on $S'$ via this isomorphism by
\[\omega_{[p,\sigma]}([p,v_{1}],[p,v_{2}]):=\omega^{S}_{\sigma}(v_{1},v_{2}),\]
for $p\in P$, $\sigma\in S$ and $v_{1},v_{2}\in T_{\sigma}S$.
This well defined, since $\omega$ is $H$-invariant by assumption. One easily sees that this indeed gives a
non-degenerate $\bar{\pi}$-2-form on $P\times_{H}S$ and that
\[d^{\bar{\pi}}\omega=0.\]

Consider the action of the gauge groupoid $G=P\times_{H}P$ on $\bar{\pi}:P\times_{H}S\to M$ given by
\[[p,q][q,\sigma]:=[p,\sigma],\]
where we remark that if $t[p,q]=\pi([q',\sigma'])$, then one can always find a representative of the class 
$[q',\sigma']$ as above. Denote this action by $\gamma$.
\begin{proposition}\label{princ2}
The action of the gauge groupoid $G=P\times_{H}P$ on $S'=P\times_{H}S$ is internally Hamiltonian.
\end{proposition}
\begin{proof}
Note that the Lie algebroid associated to $G$ is isomorphic to $(T P)/H$ (cf.\ \cite{La},\cite{Mac}). Hence the dual of the kernel 
of the anchor is isomorphic to $T^{*,\pi}P/H$, which in turn is isomorphic to $P\times_{H}\mathfrak{h}^{*}$ using 
the map 
\[(\bar{\alpha}')^{*}:T^{*,\pi}P/H\to P\times_{H}\mathfrak{h}^{*}\]
induced by the infinitesimal action $\alpha^{*}$ of $\mathfrak{h}$ on $T^{*}P$.

We give the momentum map via this isomorphism as a map
\[P\times_{H}S\to P\times_{H}\mathfrak{h}^{*}\]
defined by
\[\bar{\mu}[p,\sigma]:=[p,\mu(\sigma)].\]
This is indeed well-defined, since $\mu$ is by assumption $H$-equivariant, hence
\begin{align*}
\bar{\mu}[h\, p,h\, \sigma]&:=[h\, p,\mu(h\,\sigma)]\\
&=[h\, p,\Ad^{*}(h)\mu(\sigma)]\\
&=[p,\mu(\sigma)].
\end{align*}
Dually to $T^{\bar{\pi}}(P\times_{H}S)\simeq P\times_{H}T S$ we have an isomorphism
\[k:T^{*,\bar{\pi}}(P\times_{H}S)\to P\times_{H}T^{*}S.\]

Finally, we check that for all $X\in\Gamma(P\times_{H}\mathfrak{h})$
\begin{align*}
d^{J}\la\bar{\mu},X\ra&=d^{J}[p,\la\mu,X\ra]\\
&\mapsto[p,d^{S}\la\mu,X\ra]\\
&=[p,-\beta(X)\inner\omega^{S}]\\
&\mapsto-\gamma(X)\inner\omega,
\end{align*}
where the arrow on the second line refers to the isomorphism $k$ and in the last line we again identify
the action of $T^{\bar{\pi}}P/H$ on $P\times_{H}S$ with the action of $P\times_{H}\mathfrak{h}$ on $P\times_{H}S$
through the isomorphism $T^{\bar{\pi}}P/H\to P\times_{H}\mathfrak{h}$.

Finally we have to check equivariance of the momentum map $\bar{\mu}$. This is immediate if we again identify
$T^{\bar{\pi}}P/H$ with $P\times_{H}\mathfrak{h}$ and $T^{*,\bar{\pi}}P/H$ with $P\times_{H}\mathfrak{h}^{*}$.
\end{proof}
\end{example}

\begin{example}\label{UE}
Suppose  $\pi: E\to M$ is a  smooth complex vector bundle endowed with a Hermitian metric $h:E\fiber{\pi}{\pi}E\to\CC$.
Let $U(E)$ be the groupoid of unitary maps on the fibers \mbox{$\{E_{m}\to E_{n}\}_{m,n\in M}$}. It has a smooth structure induced from
the smooth structure on $E$ and the smooth structure on $U(n)$ (cf.\ \cite{Mac}).
There exists a smooth family of symplectic structures $\omega\in\Omega_{\pi}^{2}(E)$, given by the imaginary part
of $h$, after 
identifying $T^{\pi}E$ with $\pi^{*}E$.
\begin{proposition}
The natural action of $U(E)$ on $(E,\omega)$ is internally Hamiltonian.
\end{proposition}
\begin{proof}
Let $F_{U}(E)\subset\Hom_{M}(M\times \CC^{n},E)$ be the unitary frame bundle of $E$, i.e. the principal bundle of unitary maps of the trivial
bundle $M\times \CC^{n}$ to $E$, where $n$ is the rank of $E$. It is well-known that
\[F_{U}(E)\times_{U(n)}\CC^{n}\simeq E,\]
given by the map $(\Psi,z)\mapsto \Psi(z)$.
Moreover, one easily checks that the map 
\[F_{U}(E)\times_{U(n)}F_{U}(E)\to U(E)\]
given by $[p,q]\mapsto ([q,z]\mapsto[p,z])$ is an isomorphism of the gauge groupoid of $F_{U}(E)$ with $U(E)$.
Hence,
\[I(U(E))\simeq (F_{U}(E)\fiber{\pi}{\pi}F_{U}(E))/U(n).\]

Suppose $\omega'$ is the imaginary part of a Hermitian inner product on $\CC^{n}$. The natural action 
of $U(n)$ on $(\CC(n),\omega')$ is known to be Hamiltonian (cf.\ for example \cite{HH}).
So the proposition follows from Proposition \ref{princ2}, 
where the silent assumption was that $\omega$ is induced from
$\omega'$, as in the previous example.
\end{proof}
\end{example}

\begin{example}\label{mainexample}
Suppose a regular Lie groupoid $G$ acts on a surjective submersion $J:N\to M$. Denote the action by $\alpha:G\fiber{s}{J}N\to N$.
Denote the composition of $J:N\to M$ and the quotient map $M\to M/G$ by $\tilde{J}:N\to M/G$.
Let $p:T^{*,J}N\to N$ denote the projection and $\bar{J}:=J\circ p$. 
Note that one has a commuting diagram
\[\xymatrix{T^{*,J}N\ar[ddr]\ar[dr]^{\bar{J}}\ar[r]^{p}&N\ar[d]_{J}\ar@/^1pc/[dd]^{\tilde{J}}\\&M\ar[d]\\&M/G.}\]

There exists an induced action of the groupoid on the map $\bar{J}$ given by
\[\tilde{\alpha}(g)\eta=T^{*,J}\alpha(g)^{-1}\eta,\]
where $\eta\in (T^{*,J}N)_{s(g)}:=\bar{J}^{-1}(s(g))$. Moreover, there exists a canonical 1-form on $T^{*,J}N$ defined by 
\[\tau:=T^{*,J}p:T^{*,J}N\to T^{*,\bar{J}}(T^{*,J}N),\]
by abuse of notation ($T^{*}$ is not a functor in general.)
This gives rise to a family of symplectic forms 
\[\omega:=d^{\bar{J}}\tau\in\Omega^{\bar{J}}(T^{*,J}N).\]
\begin{proposition}
The action of $G\rra M$ on $(\bar{J}:T^{*,J}N\to M,\omega)$ is internally Hamiltonian.
\end{proposition}
\begin{proof}
We define an internal momentum map $\mu:T^{*,J}N\to\bar{J}^{*}\ker(\rho)$ by
\[\mu:=-\tilde{\alpha}^{*}\tau,\]
where we use the same notation $\tilde{\alpha}$ for the induced action of the Lie algebroid of $G\rra M$.
The fact that the action is weakly internally Hamiltonian follows from
\[d^{\bar{J}}\la\mu,\bar{J}^{*}X\ra=-d^{\bar{J}}\la\tau,\tilde{\alpha}(\bar{J}^{*}X)\ra=-\tilde{\alpha}(X)\inner d^{\bar{J}}\tau
=-\tilde{\alpha}(X)\inner\omega\]
for all $X\in\ker(\rho)$.
Equivariance of the momentum map follows from
\begin{eqnarray*}
\la\mu,\bar{J}^{*}X\ra(\tilde{\alpha}(g)\eta)&=&-\la T^{*,J}\alpha(g)^{-1}\eta,\alpha(X)_{\alpha(g)^{-1}p(\eta)}\ra\\
&=&-\la\eta,\alpha(\Ad(g)X)\ra\\
&=&\la\mu,\bar{J}^{*}(\Ad(g)X)\ra,
\end{eqnarray*}
for all $g\in G$, $\eta\in T^{*,J}N_{s(g)}$ and $X\in\ker(\rho)$.
\end{proof}
\end{example}

\begin{example}
As a corollary of the previous example, every Lie groupoid $G$ over $M$ has three canonical 
internally Hamiltonian actions associated to it.
Firstly, one has the action $G$ on the base space $M$ given by $g\cdot s(g)=t(g)$, with zero symplectic structure at each point $m\in M$.
Secondly, consider the the action on $T^{*,t}G\to M$, induced from the left action of $G$ on $t:G\to M$ by left multiplication.
Thirdly, one has the action on $T^{*,p}I_{G}$ induced from the conjugation action on $p:I_{G}\to M$. This last one is,
of course, related the coadjoint action of $G$ on $\ker(\rho)^{*}\subset T^{*,p}I_{G}$, which is internally Hamiltonian.
\end{example}

\subsection{Hamiltonian actions and momentum maps.}\label{Hamacmommap}
In this section we introduce Hamiltonian actions of Lie algebroids. A large part of the section will be devoted to
examples justifying our terminology.

Let $\AA$ be a Lie algebroid over $M$ with anchor $\rho$.
A \textbf{smooth $n$-cochain on} $\AA$ is a $C^{\infty}(M)$-multilinear antisymmetric map
\[\mu:\Gamma^{\infty}(\AA)\times\ldots\times\Gamma^{\infty}(\AA)\to C^{\infty}(M).\]
The space of smooth $n$-cochains is denoted $C^{n}(\AA)$. It is turned into a cochain complex by 
\begin{align*}d_{\AA}\mu(X_{1},\ldots,X_{n+1})&=
\sum_{i<j}(-1)^{i+j+1}\mu([X_{i},X_{j}],X_{1},\dots,\hat{X_{i}},\ldots,\hat{X_{j}},\ldots,X_{n+1})\\
&\quad+\sum_{i=1}^{n+1}(-1)^{i}\rho(X_{i})\cdot\mu(X_{1},\ldots,\hat{X_{i}},\ldots,X_{n+1}).
\end{align*}
The cohomology of the cochain complex is denoted $H^{*}(\AA)$ and is called the \textbf{Lie algebroid cohomology} of $\AA$.

\begin{remark}
This not the only type of cohomology one could associate to Lie algebroids, see e.g., \cite{CM}. The cohomology
groups discussed here are also called the (generalized) de Rham cohomology of $\AA$. One could also define
de Rham cohomology with coefficients in a Lie algebroid representation, but this is not needed in this paper.
\end{remark}

A map of Lie algebroids $\Phi:\AA\to\AA'$ induces a cochain map $\Phi^{*}:C^{*}(\AA')\to C^{*}(\AA)$ and hence a map 
$\Phi^{*}:H^{*}(\AA')\to H^{*}(\AA)$ on cohomology. 

\begin{example}
Suppose  $\AA=T M$, the Lie algebroid of the pair groupoid $M\times M$. Then $H^{*}(A)$ is the de Rham cohomology $H^{*}_{dR}(M)$ of $M$.
If $\AA=\mathfrak{g}$ is a Lie algebra, then $H^{*}(\AA)$ is the
Chevalley-Eilenberg cohomology $H^{*}(\mathfrak{g})$.
\end{example}

Suppose a regular Lie algebroid $\AA$ acts on smooth surjective submersion $J:S\to M$.
Suppose a $J$-presymplectic 2-form
\[\tilde{\omega}\in\Omega^{2}_{\tilde{J}}(S)\] 
is given. Since the action $\alpha$ is a morphism of Lie algebroids $\AA\ltimes J\to T^{\tilde{J}}S$, one has
\[d_{\AA\ltimes J}\alpha^{*}\tilde{\omega}=\alpha^{*}d^{\tilde{J}}\tilde{\omega}=0,\]
where $\AA\ltimes J$ is the action algebroid associated to the action of $\AA$ on $J:S\to M$.
Suppose the action is $J$-presymplectic. 
Using Cartan's homotopy formula this implies that
\[d^{J}(\alpha(X)\inner\,\tilde{\omega})|_{T^{J}S}=
(\alpha(X)\inner\,d^{\tilde{J}}\tilde{\omega})|_{T^{J}S}=0.\]
\begin{definition}
An action $\alpha$ of a Lie algebroid $(\pi:\AA\to M,\rho)$ on a map $J:S\to M$ is 
\textbf{Hamiltonian} if 
there exists a smooth section $\tilde{\mu}\in\Gamma^{\infty}((\AA\ltimes J)^{*})$, satisfying 
\begin{align}
d_{\AA\ltimes J}\tilde{\mu}&=-\alpha^{*}\tilde{\omega},\label{preqcond}\\
d^{J}\la\tilde{\mu},J^{*}X\ra&=-\left.(\alpha(X)\inner\,\tilde{\omega})\right|_{T^{J}S}\label{quancond}
\mbox{\textup{ for all }}X\in\Gamma^{\infty}(\AA).
\end{align}
$\tilde{\mu}$ is called a \textbf{momentum map} for the action.
\end{definition}
\begin{remark}Condition (\ref{preqcond}) is called the prequantization condition and has to be satisfied for an 
action to be prequantizable. Condition (\ref{quancond}) is called the quantization condition and has to be satisfied 
for the prequantization to be quantizable. This terminology will be justified in section 2 and 3 below. 
\end{remark}
\begin{remark}
The prequantization and quantization 
conditions state that $\tilde{\mu}\in\Gamma^{\infty}(\AA\ltimes J)$ should be the simultaneous solution of an integration problem (\ref{preqcond}) for which 
\mbox{$\alpha^{*}\tilde{\omega}\in H^{2}(\AA\ltimes J)$} is the obstruction and a lifting problem (\ref{quancond}) (which is an integration problem for each $X$)
for which the map $X\mapsto \left.(\alpha(X)\inner\,\tilde{\omega})\right|_{T^{J}S}\in H^{2}_{J}(S)$ 
forms the obstruction. In particular, if these cohomology groups are zero, then all $J$-presymplectic actions of $\AA$ on $J:S\to M$
are Hamiltonian. There exist some vanishing results for Lie algebroid cohomology (cf.\ \cite{Cra2}).
\end{remark}
\begin{lemma}
If an action of $\AA$ on $(J:S\to M,\tilde{\omega})$ is Hamiltonian, then it is internally strongly Hamiltonian.
\end{lemma}
\begin{proof}
Note that condition (\ref{quancond}) implies that that the action is internally weakly Hamiltonian, with internal momentum map
$\mu:=i^{*}\circ\tilde{\mu}$, where $i:\ker(\rho)\to\AA$ is the inclusion.

We compute the left hand side of the prequantization condition (\ref{preqcond}):
\begin{align}
d_{\AA(G)\ltimes J}\tilde{\mu}(X, Y)&=
\la\tilde{\mu},[X,Y]\ra-\alpha(X)\la\tilde{\mu},J^{*}Y\ra+\alpha(Y)\la\tilde{\mu},J^{*}X\ra.
\label{precon1}
\end{align}
If $X\in\Gamma^{\infty}(\AA)$ and $Y\in\Gamma^{\infty}(\ker(\rho))$,
then inserting (\ref{quancond}) in (\ref{precon1}) we obtain
\begin{equation}\la\tilde{\mu},[X,Y]\ra-\alpha(X)\la\tilde{\mu},J^{*}Y\ra=0.\label{eqn8}\end{equation}
But
\[\la\tilde{\mu},[X,Y]\ra=\la\tilde{\mu},\ad(X)Y\ra=\la\ad^{*}(X)\tilde{\mu},Y\ra,\]
hence (\ref{eqn8}) expresses $\AA$-equivariance of $\tilde{\mu}$.
Hence the momentum map is a lift of an internal momentum map $\mu:S\to J^{*}\ker(\rho)^{*}$, i.e.\
the diagram
\[\xymatrix{&J^{*}\AA^{*}\ar[d]^-{i^{*}}\\S\ar[r]_-{\mu}\ar[ur]^-{\tilde{\mu}}&J^{*}\ker(\rho)^{*}}\]
commutes.
\end{proof}
Suppose $\AA$ acts on $(J:S\to M,\tilde{\omega})$ in a Hamiltonian fashion.
Then we can form a perturbation of $\tilde{\omega}$ by certain exact forms without changing the fact that the 
action is Hamiltonian.
\begin{lemma}
For any $\beta\in\Omega^{1}_{\tilde{J}}(S)$. 
the action of $\AA$ on $(J:S\to M,\tilde{\omega}')$, with
\[\tilde{\omega}':=\tilde{\omega}+d^{\tilde{J}}\beta,\]
is Hamiltonian iff $\mathcal{L}_{\alpha(X)}\beta$ annihilates 
$T^{J}S$ for all $X\in\Gamma^{\infty}(\AA)$.
\end{lemma}
\begin{proof}
Define a momentum map by
\[\tilde{\mu}':=\tilde{\mu}+\alpha^{*}\beta.\]
One computes that Condition \eqref{preqcond} is satisfied:
\begin{align*}
d_{\AA\ltimes J}\tilde{\mu}'&=d_{\AA\ltimes J}\tilde{\mu}+d_{\AA\ltimes J}\alpha^{*}\beta\\
&=-\alpha^{*}\tilde{\omega}-\alpha^{*}d^{\tilde{J}}\beta\\
&=-\alpha^{*}\tilde{\omega}'.
\end{align*}
Idem dito for Condition \eqref{quancond} 
\begin{align*}
d^{J}\la\tilde{\mu}',X\ra&=d^{J}\la\tilde{\mu},X\ra+d^{J}\la\alpha^{*}\beta,X\ra\\
&=-(\alpha(X)\inner\tilde{\omega})|_{T^{J}S}-(\alpha(X)\inner d^{\tilde{J}}\beta)|_{T^{J}S}+
\mathcal{L}_{\alpha(X)}\beta|_{T^{J}S}\\
&=-(\alpha(X)\inner\tilde{\omega}')|_{T^{J}S}.
\end{align*}
\end{proof}
Moreover, given a Hamiltonian action $\alpha$ of $\AA$ on $(J:S\to M,\tilde{\omega})$ one can add certain 
closed forms to the momentum map and it is still a momentum map. Hence, the chosen momentum map is not unique.
\begin{lemma}
Suppose $\tilde{\mu}$ is a momentum map and $\beta\in\Omega_{\tilde{J}}^{1}(S)$, then
\[\tilde{\mu}':=\tilde{\mu}+\alpha^{*}\beta\]
is a momentum map for the action too
iff $d^{\tilde{J}}\beta=0$ and $\mathcal{L}_{\alpha(X)}\beta$ annihilates $T_{J}S$ for all $X\in\Gamma(\AA)$.
\end{lemma}
\begin{proof}
The proof is a calculation similar to the proof of the previous lemma.
\end{proof}
\begin{definition}
An action of a regular Lie groupoid $G$ on a smooth family of symplectic manifolds $J:S\to M$ is \textbf{Hamiltonian} 
if the induced action of the associated Lie algebroid $\AA(G)$ is Hamiltonian.
\end{definition}
\begin{example}
For the case of an action of a smooth family of Lie algebras on a smooth family of symplectic manifolds, every internal momentum map is a momentum map. In particular, strongly Hamiltonian Lie algebra actions are Hamiltonian in our terminology too.
\end{example}

\begin{example}\label{foli}
Suppose $M$ is a smooth manifold and consider an integrable distribution $T \mathcal{F}\subset T M$ as a Lie algebroid over $M$.
The differential in the Lie algebroid de Rham complex is the $\FF$-de Rham partial differential $d^{\FF}$. Consider the action of $T \mathcal{F}$ on $M\to M$.
Suppose that there exists a $\mathcal{F}$-partially closed 2-form $\tilde{\omega}$. This is trivially a $J$-presymplectic form.
A smooth section $\tilde{\mu}:S=M\to \AA^{*}(G)=T^{*}\FF$ is a momentum map iff $d^{\mathcal{F}}\tilde{\mu}=-\tilde{\omega}$ on $M$.
\end{example}

\begin{example}
Suppose a regular Lie algebroid $(\AA_{i}\to M,\rho_{i})$ acts on $(J_{i}:S_{i}\to M_{i},\tilde{\omega_{i}})$ in a Hamiltonian fashion, with momentum map $\tilde{\mu}_{i}:S\to\AA^{*}(G_{i})$ for $i=1,2$. Then the Cartesian product of algebroids $\AA_{1}\times \AA_{2}$ acts Hamiltonianly 
on 
\[(J_{1}\times J_{2}:S_{1}\times S_{2}\to M_{1}\times M_{2},\tilde{\omega}_{1}\times \tilde{\omega}_{2})\] 
with momentum map 
\[\tilde{\mu}_{1}\times\tilde{\mu}_{2}:S_{1}\times S_{2}\to\AA^{*}_{1}\times\AA^{*}_{2}.\]

As a particular example of this,
suppose $M$ is a manifold endowed with a closed 2-form $\omega$ and $\mathfrak{g}$ is a Lie algebra. Suppose $\mathcal{O}\subset\mathfrak{g}^{*}$ is a coadjoint orbit. Consider the trivial Lie algebroid over 
$M$ with fiber $\mathfrak{g}$. It is the Cartesian 
product $T M\times \mathfrak{g}$ of $T M$ as a Lie algebroid over $M$ and $\mathfrak{g}$ as a Lie algebroid over a point. The inclusion of the
coadjoint orbit $i:\mathcal{O}\hookrightarrow\mathfrak{g}^{*}$ is a momentum map for the coadjoint action of $\mathfrak{g}$ on $\mathcal{O}$.
If $\beta\in\Omega^{1}(M)$ satisfies $d\beta=-\omega$, then $\beta:M\to T^{*}M$ is a momentum map for the action
of $T M$ on $(M\to M,\omega)$.
Hence 
\[\beta\times i:M\times\mathcal{O}\to T^{*}M\times \mathfrak{g}^{*}\] 
is a momentum map for the action of $T M\times\mathfrak{g}$ on $pr_{1}:M\times\mathcal{O}\to M$. 
\end{example}

\begin{example}\label{princ3}
Consider the situation of Example \ref{princ}, in which a Lie group $H$ acts in a Hamiltonian fashion on a symplectic manifold 
$(S,\omega^{S})$,
with momentum map $\mu:S\to\mathfrak{h}^{*}$. Suppose $\pi:P\to M$ is a principal $H$-bundle. Proposition \ref{princ2} 
states that the action of the gauge groupoid $G=P\times_{H}P\rra M$ on $J:S':=P\times_{H}S\to M$ is internally Hamiltonian. Given a 
connection on $P$ we shall extend the symplectic form $\omega\in\Omega_{J}(S')$ on $S'$ to a $J$-presymplectic 
form $\tilde{\omega}\in\Omega^{2}_{\tilde{J}}(S)$. Then we shall see that the action 
is Hamiltonian with respect to a well-chosen momentum map.

Suppose $\tau\in\Gamma^{\infty}(\Wedge^{1}(P)\otimes\mathfrak{h})$ is a Lie algebra-valued connection 1-form on $P$.
After identifying $T(P\times_{H}S)\simeq(T P\times_{H}T S)/\mathord\sim$ as in Example \ref{princ}, define
\begin{align*}\tilde{\omega}_{[p,\sigma]}\bigl([w_{1},v_{1}],[w_{2},v_{2}]\bigr)&:=
\omega^{S}_{\sigma}\bigl((v_{1}-\beta(\tau(w_{1})),v_{2}-\beta(\tau(w_{2}))\bigr)
-\la\mu_{\sigma}, F_{p}(w_{1},w_{2})\ra,
\end{align*}
where $F$ is the $\mathfrak{h}$-valued curvature 2-form on $P$.
\begin{lemma}\label{wdom}
$\tilde{\omega}$ is a well-defined 2-form in $\Omega_{\tilde{J}}^{2}(S')$.
\end{lemma}
\begin{proof}
Note that for $X,Y\in\mathfrak{h}$
\[\tilde{\omega}\bigl([\alpha(X),\beta(X)],[\alpha(Y),\beta(Y)]\bigr)=0.\]
Moreover, for all $h\in H$ 
\[\tilde{\omega}_{[h\,p,h\,\sigma]}\bigl([h\,w_{1},h\,v_{1}],[h\,w_{2},h\,v_{2}]\bigr)=
\tilde{\omega}_{[p,\sigma]}\bigl([w_{1},v_{1}],[w_{2},v_{2}]\bigr),\]
since $\omega^{S}$ is $H$-invariant and
\begin{align*}\la\mu_{h\,\sigma}, F_{h\,p}(h\,w_{1},h\,w_{2})\ra
&=\la\Ad^{*}(h)\,\mu_{\sigma},\Ad(h)\,F_{p}(w_{1},w_{2})\ra\\
&=\la\mu_{\sigma},F_{p}(w_{1},w_{2})\ra,\end{align*}
by $H$-equivariance of $\tau$ and $\mu$.
\end{proof}
We shall omit here the proof that $\tilde{\omega}$ is closed since we shall later see that it is the curvature 2-form of 
a connection on a line bundle over $P\times_{H}S$. Obviously, $\tilde{\omega}$ restricts to $\omega$ on the vertical tangent space.

\begin{proposition}\label{princ6}
The action of the gauge groupoid $G=P\times_{H}P\rra M$ on \mbox{$(J:S'\to M,\tilde{\omega})$} is Hamiltonian.
\end{proposition}
\begin{proof}
Define a momentum map $\tilde{\mu}:P\times_{H}S\to \AA^{*}\simeq T^{*}P/H$ by
\[\tilde{\mu}:=\la\mu,\tau\ra,\]
where $\mu$ is the momentum map for the action of $H$ on $(S,\omega^{S})$.
This is well-defined, since
\begin{align*}\la\mu,\tau\ra(h\,p,h\,\sigma)&=\la \Ad^{*}(h)\mu(\sigma),\Ad(h)\tau(p)\ra\\
&=\la\mu,\tau\ra(p,\sigma).\end{align*}
One easily sees that $\tilde{\mu}$ restricts to $\bar{\mu}$ (cf.\ Example \ref{princ}) on the vertical tangent space.
For $H$-equivariant vector fields $w_{1},w_{2}$ on $P$, we compute
\begin{align*}
d_{P}\tilde{\mu}(w_{1},w_{2})&=\la\mu,d_{P}\tau(w_{1},w_{2})\ra\\
&=\la\mu, F(w_{1},w_{2})\ra-\la\mu,[\tau(w_{1}),\tau(w_{2})]_{\mathfrak{h}}\ra\\
&=\la\mu, F(w_{1},w_{2})\ra+\beta(\tau(w_{1}))\cdot\la\mu,\tau(w_{2})\ra\\
&=\la\mu, F(w_{1},w_{2})\ra-\omega(\tau(w_{1}),\tau(w_{2})),
\end{align*}
where the second equality follows from the curvature formula
\[F=d\tau+[\tau,\tau]_{\mathfrak{h}},\]
the third equality follows from $H$-equivariance of $\mu$ and
the last equality follows from the fact that $\mu$ is a momentum map.
Hence one has $d_{A\ltimes J}\tilde{\mu}=-\gamma^{*}\tilde{\omega}$, where $\gamma$ denotes the action of $P\times_{H}P$
on $P\times_{H}S$.

We check the quantization condition (\ref{quancond}) for $\tilde{\mu}$. Identify 
\[T^{J}(P\times_{H}S)\simeq P\times_{H}T S,\] 
as in Example \ref{princ}. For $w\in\mathfrak{X}(P)$ we compute
\begin{align*}
d^{S}\la\tilde{\mu},w\ra&=d^{S}\la\mu,\tau(w)\ra\\
&=-\beta(\tau(w))\inner\omega^{S},
\end{align*}
from which we conclude that
\[d^{J}\la\tilde{\mu},w\ra=-\gamma(w)\inner\tilde{\omega}.\]
\end{proof}
\end{example}

\begin{example}
Suppose $\pi:E\to M$ is a complex vector bundle with Hermitian structure $h$.
Consider the action of $U(E)$ on $E$ as in Example \ref{UE}. Let $F_{U}(E)\to M$ be the frame bundle of $E$.
Suppose $F_{U}(E)$ is endowed with a connection $\tau$. Then we can extend $\omega$ (cf.\ Example \ref{UE}) to a closed form
$\tilde{\omega}$ on $E$ as in the above Example \ref{princ3}. As a consequence of Proposition \ref{princ6} we have 
the following
\begin{corollary}
The action of $U(E)$ on $(E,\tilde{\omega})$ is Hamiltonian.
\end{corollary}
\end{example}
\begin{remark}
Note that the input for a Hamiltonian Lie groupoid action is not just a smooth bundle of coadjoint orbits.
One can start with a smooth bundle of coadjoint orbits $J:S=\{\mathcal{O}_{m G}\}_{m G\in M/G}\to M$ 
on which one has an internal Hamiltonian action (cf.\ Example \ref{coadj}), but the difficulty is to find a 
suitable $J$-presymplectic form 
on it, such that the coadjoint action is Hamiltonian. For this reason, in the context of Lie groupoids 
there is not a perfect orbit method or philosophy \`{a} la Kirillov. We shall further discuss this 
in Section \ref{orbitmethod}.
\end{remark}

\section{Prequantization of Hamiltonian Lie algebroid actions}
\subsection{Representations of Lie algebroids on vector bundles.}
In this section we introduce the notion of a representation of a Lie algebroid with a 
base manifold possibly different from the base manifold of $\AA$. Finally we introduce a certain
Picard group associated to Lie algebroid representations on line bundles and show that there exists an exact sequence
involving this Picard group.

Suppose $(p:\AA\to M,\rho)$ is a Lie algebroid
and $\alpha:\Gamma^{\infty}(\AA)\to\XX^{\infty}_{\tilde{J}}(S)$ is an action of $\AA$ on $J:S\to M$. 
Suppose $E\to S$ is a smooth vector bundle over $S$. 
Let $\mathcal{D}(E)$ be the Lie algebroid whose sections are first-order differential operators $D$ on 
$E\to S$ that for all $f\in C^{\infty}_{\tilde{J}}(S)$ and $\sigma\in\Gamma^{\infty}(E)$ satisfy
\[D(f \sigma)=f D\sigma+\Theta(D)f\sigma\]
for a map $\Theta:\mathcal{D}(E)\to T^{\tilde{J}}S$, which is the anchor (see \cite{Mac}).

\begin{definition}
A \textbf{$\AA$-connection on a complex vector bundle} $p:E\to S$,
is a map of vector bundles $\pi:\AA\ltimes J\to\mathcal{D}(E)$, such that $\alpha=\Theta\circ \pi$.
If $\pi$ preserves the Lie bracket, then it is called a \textbf{representation} or \textbf{flat $\AA$-connection}.
\end{definition}
\begin{remark}
This is a more general notion of a representation than usual, in the sense that we allow a base manifold which is
not $M$.
We shall see that prequantization defines a representation in this way. Quantization gives a representation
in the narrow sense, i.e. on a smooth vector bundle over $M$.
\end{remark}

\begin{definition}
A representation $\pi:\AA\ltimes J\to\mathcal{D}(E)$ is \textbf{Hermitian} with respect to a Hermitian metric $g$ on $E\to M$ if
\[g(\pi(X)\sigma,\tau)+g(\sigma,\pi(X)\tau)=\alpha(X)g(\sigma,\tau).\]
\end{definition}

\begin{proposition}
Any Hermitian representation of a Lie algebroid $\AA$ on a line bundle $L\to S$ is of the form 
\[\pi(J^{*}X):=\nabla_{\alpha(X)}-2\pi i\la \tilde{\mu},J^{*}X\ra,\]
where $\alpha$ is the action of $\AA$ on $J:S\to M$, $\tilde{\mu}\in\Gamma^{\infty}((\AA\ltimes J)^{*})$
and $\nabla$ a Hermitian $T^{\tilde{J}}S$-connection on $S$.
\end{proposition}
This proposition follows from the following lemma.
\begin{lemma}
Let $(p:\AA'\to M,\rho)$ be a Lie algebroid and let $\nabla$ be a Hermitian $\AA$-connection on a smooth complex vector bundle $E\to M$.
Then there exist a Hermitian connection $\nabla^{E}$ on $E\to M$ and a section $\mu$ of $\AA^{*}\otimes\End(E)$ such that $\nabla$ is of the form
\[\nabla_{X}=\nabla^{E}_{X}-2\pi i\langle\mu,X\rangle,\]
for all $X\in\Gamma^{\infty}(\AA)$.
\end{lemma}
\begin{proof}
Let $\nabla^{E}$ be any Hermitian connection on $E$. It is well known that such a connection always exists. Consider the associated 
$\AA$-connection defined by $\tilde{\nabla}:X\mapsto\nabla_{\rho(X)}^{E}$. Now, 
\begin{align*}
(\nabla_{X}-\tilde{\nabla}_{X})(f s)&=f\nabla_{X}s+\rho(X)f s-f\nabla^{E}_{\rho(X)}s-\rho(X)f s\\
&=f(\nabla_{X}-\tilde{\nabla}_{X})s,
\end{align*}
hence $\nabla_{X}-\tilde{\nabla}_{X}$ is a zeroth order differential operator on $E$, i.e. $\nabla_{X}-\tilde{\nabla}_{X}\in\End(E)$.
Moreover, $\nabla-\tilde{\nabla}$ is $C^{\infty}(M)$-linear, in the sense that 
\[\nabla_{f X}-\tilde{\nabla}_{f X}=f(\nabla_{X}-\tilde{\nabla}_{X}),\]
Thus $\nabla-\tilde{\nabla}\in \Gamma^{\infty}(\AA^{*}\otimes\End(E))$ by the Serre-Swan theorem.
\end{proof}

Isomorphism classes of smooth complex line bundles on a manifold $S$ form a 
group $\Pic_{\tilde{J}}(S)$ under the tensor 
product, with the trivial rank one line bundle as a unit, and inverse $[L]^{-1}=[L^{*}]:=[\Hom_{S}(L,\CC)]$. If $\AA$ acts on a map $J:S\to M$, then 
one can extend this structure to the set of isomorphism classes of Hermitian $\AA$-representations 
on smooth complex line bundles over $S$. The product of $\pi:\AA\ltimes J\to \mathcal{D}(L)$ and $\pi':\AA\ltimes J\to\mathcal{D}(L')$ is defined by
\[\pi\otimes 1+1\otimes\pi':A\ltimes J\to \mathcal{D}(L\otimes L'),\]
and the inverse is given by
\[\la\pi^{-1}(X)s^{*},s\ra=-\langle s^{*},\pi(X)(s)\rangle+d_{\AA}\langle s^{*},s\rangle.\]
The Hermitian structure on the tensor product $L_{1}\otimes L_{2}$ of two line bundles $(L_{1},g_{1})$, $(L_{2},g_{2})$ is given by the 
formula
\[g(v_{1}\otimes w_{1},v_{2}\otimes w_{2})=g_{1}(v_{1},v_{2})g_{2}(w_{1},w_{2}).\]
\begin{definition}
The \textbf{Hermitian Picard group $\Pic_{\AA}(J)$ of the Lie algebroid action} of $\AA$ on $J:S\to M$ is the group of isomorphism classes of 
Hermitian representations of $\AA$ on line bundles over $S$, with product and inverse as described above.
\end{definition}
\begin{proposition}
There is an exact sequence of groups
\begin{align}\label{ses}
0\to H^{1}(\AA\ltimes J)\to\Pic_{\AA}(J)\to\Pic(S)\overset{c_{1}^{\AA\ltimes J}}\To H^{2}(\AA\ltimes J).
\end{align}
\end{proposition}
\begin{proof}
The second arrow sends a closed section $\mu\in\Gamma^{\infty}((\AA\ltimes J)^{*})$ to the representation 
$X\mapsto \alpha(X)-2\pi i\langle\mu,X\rangle$ on the trivial line bundle.
This is well-defined and injective: suppose $\mu,\mu'\in\Gamma^{\infty}((\AA\ltimes J)^{*})$ give rise to isomorphic representations, i.e. 
there is a $f\in C^{\infty}(S)$ such that
\[(\alpha(X)-2\pi i\langle\mu',J^{*}X\rangle)(f \sigma)=f (\alpha(X)-2\pi i\langle\mu,J^{*}X\rangle) \sigma,\] 
for all $X\in \Gamma^{\infty}(\AA)$ and all $\sigma\in C^{\infty}_{\tilde{J}}(S)=\Gamma^{\infty}(L)$. Using the Leibniz rule and 
$\langle d_{\AA\ltimes J}(f),X\rangle=\alpha(X)f$ we obtain
\[\langle d_{\AA}f,J^{*}X\rangle=\langle\mu-\mu',J^{*}X\rangle.\]
for all $X\in \Gamma^{\infty}(\AA)$. Thus the two representations are isomorphic iff there exists an $f\in C^{\infty}_{\tilde{J}}(S)$ such that 
$d_{A}f=\mu-\mu'$. The third arrow forgets the representation, so the sequence is exact at $\Pic_{\AA}(S)$. 
The last arrow is the first $\AA$-Chern class map
$c_{1}^{\AA\ltimes J}([L]):=[\alpha^{*}K]$, where $K$ denotes the curvature 2-form of any connection $\nabla$ on $L\to M$. If it is zero in $H^{2}(\AA\ltimes J)$, then $L$ carries a Hermitian $\AA$-representation, cf. Theorem \ref{preqrep}.
\end{proof}
\begin{remark}
One can generalize the notion of (higher) Chern classes of complex vector bundles to characteristic classes for 
complex representations of Lie algebroids (cf.\ \cite{Cra2}). 
\end{remark}

\subsection{Prequantization line bundles and longitudinal \v{C}ech cohomology.}
In this section we shall discuss how the class of a $J$-presymplectic form 
$[\tilde{\omega}]\in H^{2}_{\tilde{J},dR}(S)$ determines a class in the longitudinal \v{C}ech cohomology.
This cohomology is defined analogously to the usual \v{C}ech cohomology (cf. \cite{BoTu}).
Then we give a criterion for $\tilde{\omega}$ to be the curvature of a Hermitian connection on a complex line bundle.

Suppose $\mathcal{F}$ is a regular foliation of $S$. Consider the projection map on the orbit space $\tilde{J}:S\to S/\mathcal{F}$.
Suppose $\mathcal{U}$ is a countable good foliation covering for $S$ (i.e.\ for all $U\in\mathcal{U}$ the foliation restricted
to $U$ is diffeomorphic to a contractible open subset of $\RR^{q}\times \RR^{n-q}$, where $n=\dim(M)$ and $q$ the dimension of the foliation. Let $I$ be a countable ordered 
index set for $\mathcal{U}$. Denote the intersection of $k$ sets $U_{i_{1}},\ldots, U_{i_{k}}$ by $U_{i_{1}\ldots i_{k}}$
for $i_{1},\ldots, i_{k}\in I$.
For $k\in\ZZ_{\geq 0}$ let $C_{\tilde{J}}^{k}(\mathcal{U},\RR)$ be the vector space of smooth functions on
$(k+1)$-fold intersections $U_{i_{1},\dots i_{k+1}}$ (where $i_{1}<\dots<i_{k+1}$)
which are locally constant along the leaves of the foliation by $\tilde{J}$. Define a map 
\[\delta_{k}:C_{\tilde{J}}^{k}(\mathcal{U},\RR)\to C_{\tilde{J}}^{k+1}(\mathcal{U},\RR)\]
by the usual formula
\begin{equation}\label{del}
\delta_{k}(f)|_{U_{i_{1}\ldots i_{k+1}}}:=\sum_{j=1}^{k+1}(-1)^{j}f|_{U_{i_{1}\ldots\hat{i}_{j}\ldots i_{k+1}}}.
\end{equation}
One checks that $\delta^{2}=0$. The cohomology of the complex is independent of the chosen good foliation cover and
we call it the \textbf{longitudinal \v{C}ech cohomology} and denote it by $\check{H}^{*}_{\tilde{J}}(S,\RR)$

Consider the foliation \v{C}ech-de Rham double complex defined by
\[C^{k,l}:=\prod_{i_{1}<\ldots<i_{k+1}}\Omega_{\tilde{J}}^{l}(U_{i_{1}\ldots i_{k+1}}),\]
with
\[\delta_{k,l}:C^{k,l}\to C^{k+1,l}\]
the straightforward generalization of \eqref{del} and
\[d^{\tilde{J}}_{k,l}:C^{k,l}\to C^{k,l+1}\]
the restriction of $d^{\tilde{J}}$ to the $(k+1)$-fold intersections.

The augmented double complex (ignore the fact that some arrows are dotted), partly shown here,
\[\xymatrix{0\ar[r]&\Omega^{2}_{\tilde{J}}(S)\ar@{.>}[r]&C^{0,2}\ar[r]^{\delta}&C^{1,2}\ar[r]^{\delta}&C^{2,2}\\
0\ar[r]&\Omega^{1}_{\tilde{J}}(S)\ar[u]_{d^{\tilde{J}}}\ar[r]&C^{0,1}\ar@{.>}[u]_{d^{\tilde{J}}}\ar@{.>}[r]^{\delta}&C^{1,1}
\ar[u]_{d^{\tilde{J}}}\ar[r]^{\delta}&C^{2,1}\ar[u]_{d^{\tilde{J}}}\\
0\ar[r]&\Omega^{0}_{\tilde{J}}(S)\ar[u]_{d^{\tilde{J}}}\ar[r]&C^{0,0}\ar[u]_{d^{\tilde{J}}}\ar[r]^{\delta}&C^{1,0}
\ar@{.>}[u]_{d^{\tilde{J}}}\ar@{.>}[r]^{\delta}&C^{2,0}\ar[u]_{d^{\tilde{J}}}\\
&&C^{0}_{\tilde{J}}(\mathcal{U},\RR)\ar[u]\ar[r]^{\delta}&C^{1}_{\tilde{J}}(\mathcal{U},\RR)\ar[u]\ar[r]^{\delta}&
C^{2}_{\tilde{J}}(\mathcal{U},\RR)\ar@{.>}[u]\\
&&0\ar[u]&0\ar[u]&0\ar[u]}\]
can be used to prove 
\begin{proposition}
There exists an isomorphism $H^{*}_{\tilde{J},dR}(S)\simeq\check{H}^{*}_{\tilde{J}}(S,\RR)$ between the foliation de Rham cohomology and
the longitudinal \v{C}ech cohomology. 
\end{proposition}
The proof is analogous to the proof with the usual \v{C}ech-de Rham complex (cf.\ \cite{BoTu}).

Let $[\tilde{\omega}]\in H^{2}_{\tilde{J},dR}(S)$ be the class of a $J$-presymplectic form. We shall concretely realize 
the above isomorphism to associate a degree 2 longitudinal \v{C}ech cohomology class to $[\tilde{\omega}]$.
We shall follow the dotted arrows in the above diagram.
Suppose $\mathcal{U}$ is a good foliation covering for $S$. Since $d^{\tilde{J}}\tilde{\omega}=0$, for each $U_{j}\in\mathcal{U}$ 
there exists an $\eta_{j}\in\Omega^{1}_{\tilde{J}}(U_{j})$ such that $d^{\tilde{J}}\eta_{j}=\tilde{\omega}|_{U_{j}}$.
Since for all $U_{j},U_{k}\in\mathcal{U}$ we have $d^{\tilde{J}}(\eta_{j}-\eta_{k})=0$ on the intersection 
$U_{j k}$, there exists an $f_{j k}\in C^{\infty}(U_{j k})$ such that $d^{\tilde{J}}f_{j k}=\eta_{j}-\eta_{k}$.
One easily checks that $d^{\tilde{J}}f_{j k}+d^{\tilde{J}}f_{k l}-d^{\tilde{J}}f_{j l}=0$ on 
$U_{j k l}$ and $\delta(a_{j k l})=0$. Hence 
\[a:=\{a_{j k l}:=f_{j k}+f_{k l}-f_{j l}\}_{j, k, l\in I}\]
defines a class $[a]$ in $\check{H}^{2}_{\tilde{J}}(\mathcal{U},\RR)$.

There is an obvious definition of longitudinal \v{C}ech cohomology $\check{H}^{*}_{\tilde{J}}(S,\ZZ)$ 
with values in $\ZZ$. But, one easily sees that $\check{H}^{*}_{\tilde{J}}(S,\ZZ)\simeq \check{H}^{*}(S,\ZZ)$,
since a cocycle which is integral valued and continuous is locally constant. We call a class $[\tilde{\omega}]\in H^{2}_{\tilde{J},dR}(S)$ integer if the associated class
in $\check{H}^{2}_{\tilde{J}}(S,\RR)$ is in the image of the canonical map
\[\check{H}^{2}_{\tilde{J}}(S,\ZZ)\To \check{H}^{2}_{\tilde{J}}(S,\RR).\]

\begin{theorem}\label{exline}
A $\tilde{J}$-closed form $\tilde{\omega}\in\Omega^{2}_{\tilde{J}}(S)$ is the curvature 2-form of a $\tilde{J}$-partial 
Hermitian connection on a complex line bundle $L\to S$ iff $[\tilde{\omega}]\in H^{2}_{\tilde{J},dR}(S)$ is integral. 
\end{theorem}
\begin{proof}
$(\Rightarrow)\quad$ 
Suppose a line bundle $L\to S$, Hermitian metric $h$ and a Hermitian connection $\nabla$ are given, such that $K^{\nabla}=\tilde{\omega}$.
Suppose $\{(U_{j},s_{j})\}_{j\in I}$ form a normalized trivialization of $L\to S$, in the sense that
$s_{j}:U_{j}\to L|_{U_{j}}$ is a section for all $j\in I$ such that $h(s_{j},s_{j})=1$. This gives rise to a cocycle 
$\{c_{j k}:U_{j k}\to U(1)\}_{j, k\in I}$ defined by
\[s_{k}=c_{j k}s_{j}\]
for all $j, k\in I$. 

To the curvature form $\tilde{\omega}$ of the connection $\nabla$ is associated 
a \v{C}ech class as above. The local $\tilde{J}$-forms
$\eta_{j}$ ($j\in I$) give the partial connection with curvature $\tilde{\omega}$ by the formula
\[\nabla s_{j}=2 \pi i \eta_{j}\cdot s_{j}.\]
From this formula one computes, using the Leibniz rule for connections, that
\[d^{\tilde{J}}f_{j k}=\eta_{k}-\eta_{j}=\frac{1}{2 \pi i}\frac{d^{\tilde{J}}c_{j k}}{c_{j k}}.\]
One can easily show that the fact that $\nabla$ is Hermitian implies that the function $f_{j k}$ must be 
real-valued. Hence for all $j, k\in I$ 
\[c_{j k}=e^{2\pi i ((f_{j k}+ d_{j k})},\]
for a function $d_{j k}:U_{j k}\to\RR$ locally constant along the leaves. The $d_{j k}$ constitute a \v{C}ech 
1-cocycle in $b\in C_{\tilde{J}}^{1}(\mathcal{U},\RR)$. From the fact that $c_{j k}c_{k l}c_{j l}^{-1}=1$ we deduce that
\[(f_{j k}+ d_{j k})+(f_{j k}+ d_{j k})-(f_{j k}+ d_{j k})\in\ZZ,\]
hence $a-\delta(b)\in\ZZ$, which implies that $[a]$ is integer.
\\
$(\Leftarrow)\quad$ Suppose an integer class in $H^{2}_{J,dR}(S)$ is given. There exist an associated class
in $\check{H}^{2}_{\tilde{J}}(\mathcal{U},\RR)$. Choose a representative of this class such that
the functions $a_{j k l}:=f_{j k}+f_{k l}-f_{j l}$ (as above) have integer value.
For all $j, k\in I$ define 
\[c_{j k}:=e^{2\pi i f_{j k}}\]
This defines is a cocycle, since $a_{j k l}$ is integral, which gives a smooth complex line bundle
\[L=\bigl(\bigcup_{j\in I}(U_{j}\times\CC)\bigr)/\mathord\sim,\]
where the equivalence relation is given by $U_{j}\times\CC\ni(m,z)\sim(m,c_{j k}z)\in U_{k}\times\CC$
whenever $m\in U_{j k}$. The connection is given by
\[\nabla s_{j}=2 \pi i \eta_{j}\cdot s_{j},\]
where $s_{j}$ denote section $U_{j}\times\CC$ which is constantly equal to 1. The Hermitian structure is given by
\[h((m,z_{1}),(m,z_{2}))=\bar{z}_{1}z_{2}.\]
A computation proves that $\nabla$ is Hermitian with respect to $h$.
\end{proof}

\begin{remark}
One should relate this to the well-known fact that
\[c_{1}:\Pic(S)\overset{\simeq}{\To}\check{H}^{2}(S,\ZZ)),\]
where $c_{1}([L]):=[\tilde{\omega}]$ is the first Chern class, which equals the class of the curvature 2-form of any connection on $L$.
Moreover, as we remarked before, $\check{H}^{2}(S,\ZZ)\simeq\check{H}_{\tilde{J}}^{2}(S,\ZZ)$, hence the above proof is very similar to
the proof of the fact that $c_{1}$ is an isomorphism. It is repeated here since $\tilde{\omega}$ gives rise to an element in $\check{H}_{\tilde{J}}^{2}(S,\RR)$ and not in $\check{H}^{2}(S,\RR)$ and for expositionary purposes. Summarizing one has the following commuting diagram
\[\xymatrix{0\ar[r]&\Pic_{\AA}(J)\ar[r]&\Pic(S)\ar[d]_{\simeq}\ar[rrr]^{c_{1}^{\AA}}&&&H^{2}(\AA\ltimes J)\\
&&\check{H}^{2}(S,\ZZ)\ar[r]^{\simeq}&\check{H}^{2}_{\tilde{J}}(S,\ZZ)\ar[r]
&\check{H}^{2}_{\tilde{J}}(S,\RR)\ar[r]^{\simeq}&H^{2}_{\tilde{J},dR}(S)\ar[u]^{\alpha^{*}}.}\]
\end{remark}

\subsection{Prequantization representations.}\label{preq}
In this section we prove that under suitable assumptions there exists a prequantization representation associated to a 
Hamiltonian Lie algebroid action. Next, we discuss some examples and some properties of the prequantization representation. 

Suppose a regular Lie algebroid $\AA$ over $M$ acts on a smooth map $J:S\to M$. 
Let $\alpha:\Gamma^{\infty}(\AA)\to \XX_{\tilde{J}}^{\infty}(S)$ denote the action.
Suppose that $S$ is endowed with a $J$-presymplectic $\tilde{J}$-2-form $\tilde{\omega}\in\Omega_{\tilde{J}}^{2}(S)$. 
Suppose that the action is Hamiltonian with momentum map 
\[\tilde{\mu}:S\to J^{*}\AA^{*}.\]

Suppose, furthermore, that there exists a smooth complex line bundle $L\to S$ 
with a Hermitian metric $h$
and a $\tilde{J}$-partial Hermitian connection $\nabla^{L}$, such that the curvature $J$-2-form $K^{\nabla_{L}}$ equals $\tilde{\omega}$.
The triple $(L\to S,\nabla^{L},h)$ is called a prequantization of the Hamiltonian action of $\AA$ on $(J:S\to M,\tilde{\omega})$.
We have seen in Theorem \ref{exline} that a prequantization exists if and only if the cohomology class of 
$\tilde{\omega}$ is integral.

\begin{theorem}\label{preqrep}
There exists a Hermitian representation of the Lie algebroid $\AA$ on \mbox{$(L\to M, h)$} given by
\[\pi(X):=\nabla_{\alpha(X)}-2\pi i\la \tilde{\mu},J^{*}X\ra.\]
\end{theorem}
\begin{remark}
Note that this formula is a generalization of the well-known Kostant formula in classical prequantization theory.
The fact that it also applies to Lie algebroids was also used in \cite{WeZa}.
Actually one only needs the prequantization condition (\ref{preqcond}), which is equivalent to $c_{1}^{\AA}(L)=0\in H^{2}(\AA\ltimes J)$.
From exactness of the sequence \ref{ses} the theorem follows at once. 
\end{remark}
\begin{proof}
 For $X,Y\in\Gamma^{\infty}(\AA)$ one computes
\begin{align*}
[\pi(X),\pi(Y)]_{\mathcal{D}(L)}-\pi([X,Y]_{\AA})
&=[\nabla_{\alpha(X)},\nabla_{\alpha(Y)}]\\
&\quad+2\pi i \alpha(Y)\la\tilde{\mu},J^{*}X\ra-2\pi i \alpha(X)\la\tilde{\mu},J^{*}Y\ra\\
&\quad-\nabla_{\alpha[X,Y]}+2\pi i\la\tilde{\mu},J^{*}[X,Y]\ra\\
&=2\pi i K(\alpha(X),\alpha(Y))+2 \pi i d_{A\ltimes J}\tilde{\mu}(X,Y)\\
&=2\pi i\omega(\alpha(X),\alpha(Y))+2 \pi i d_{A\ltimes J}\tilde{\mu}(X,Y)\\ 
&=0.\end{align*}
So $\pi$ is a homomorphism of Lie algebroids.

The representation being Hermitian is proven by computation. For $\sigma,\tau\in\Gamma^{\infty}(L)$ and $X\in\Gamma^{\infty}(\AA)$,
\begin{align*}
h(\pi(X)\sigma,\tau)+h(\sigma,\pi(X)\tau)&=h(\nabla_{\alpha(X)}\sigma,\tau)-h(2\pi i\la\mu,J^{*}X\ra\sigma,\tau)\\
&\quad+\,h(\sigma,\nabla_{\alpha(X)}\tau)-h(\sigma,2\pi i\la\mu,J^{*}X\ra\tau)\\
&=\alpha(X)h(\sigma,\tau)
\end{align*}
since the connection is Hermitian and the metric $h$ is sesquilinear.
\end{proof}
\begin{definition}
The above representation $(L\to M,h,\pi)\in\Pic_{\AA}(J)$ is the \textbf{prequantization representation} of the Hamiltonian action of $\AA$
on $J:S\to M$.
\end{definition}
\begin{example}\label{flatcon}
The tangent bundle $T M$ of a smooth manifold $M$ is a Lie algebroid over $M$. It trivially acts on $J=id:M\to M$.
The prequantization procedure boils down to a standard situation in differential geometry.
Suppose $M$ is endowed with an integral closed 2-form $\tilde{\omega}\in\Omega^{2}(M)$ (the $J$-presymplectic form). 
As we have seen, a momentum map 
for this action is a 1-form $\mu\in\Omega^{1}(M)$ satisfying $d\mu=-\tilde{\omega}$. A prequantum line bundle is a 
complex line bundle
$L\to M$ endowed with a Hermitian connection $\nabla$ whose curvature equals $\omega$. The prequantization representation
of $\AA=T M$ is the flat connection $\nabla-2\pi i\mu$.

For a regular integrable distribution as a Lie algebroid and its associated foliation $\FF$ (cf.\ Example \ref{foli}) a similar
reasoning holds with the differential $d$ replaced by a partial differential along the leaves of the foliation $\FF$.
\end{example}

\begin{example}\label{bundlelie2}
Suppose $p:\mathfrak{g}\to M$ is a smooth family of Lie algebras $\mathfrak{g}_{m}$ ($m\in M$) as in 
Example \ref{bundlelie}. Suppose it acts in a Hamiltonian fashion on a smooth bundle $J:S\to M$ of coadjoint orbits 
$S_{m}:=\mathcal{O}_{m}\subset\mathfrak{g}_{m}^{*}$. Then the inclusion $S\to\mathfrak{g}^{*}$ is a momentum map. 
We have a prequantization line bundle
if we can paste prequantization line bundles $L_{m}\to S_{m}$ for each $\mathfrak{g}_{m}$ into a smooth 
bundle $L\to M$. A Hermitian representation on $L\to M$ is then given by 
$X\mapsto -2 \pi i\la\tilde{\mu},J^{*}X\ra$,
where $\tilde{\mu}:S\to\mathfrak{g}^{*}$ is the inclusion.
\end{example}

\begin{example}\label{princ4}
Let a Lie group $H$ act on a symplectic manifold $(S,\omega^{S})$ in Hamiltonian
fashion, with momentum map $\mu:S\to\mathfrak{h}^{*}$, and let $\pi:P\to M$ be a principal bundle endowed with an $\mathfrak{h}$-valued connection 1-form $\tau$. 
The connection induces a decomposition of $T P$ into a direct sum $\mathcal{H}\oplus\mathcal{V}$ of a horizontal
bundle $\mathcal{H}:=\ker(\tau)$ and a vertical bundle $\mathcal{V}:=\ker(T \pi)$.
In Example \ref{princ3} we defined a $J$-presymplectic 2-form $\tilde{\omega}$ on $S':=P\times_{H}S$ 
and proved that the action of the gauge groupoid $P\times_{H}P\rra M$ on $(S',\tilde{\omega})$ is Hamiltonian.

Suppose $(\pi_{L}:L\to S,\nabla^{L},h)$ is a prequantization for the action of $H$ on $S$.
Consider the line bundle $P\times_{H}L\to P\times_{H}S$. We shall show that this forms a prequantization line bundle for
the action of $P\times_{H}P\rra M$ on $S'\to M$.

Firstly, we shall explain the line bundle structure on $\pi:P\times_{H}L\to P\times_{H}S$.  The map
\[\pi([p,z]):=[p,\pi_{L}(z)],\]
is well-defined, since $\pi_{L}$ is $H$-equivariant.
Addition is defined by finding representatives with equal first entry (this is always possible) 
and then adding the second entry, i.e.
\[[p,z_{1}]+[p,z_{2}]=[p,z_{1}+z_{2}].\]
Scalar multiplication is defined by scalar multiplication on the second entry
\[\lambda[p,z]=[p,\lambda z].\]

A section $\theta\in\Gamma^{\infty}(P\times_{H}L)$ is represented by a pair $(\theta_{1},\theta_{2})$ of $H$-equivariant maps
$\theta_{1}:P\times S\to P$ and $\theta_{2}:P\times S\to L$, such that $\theta_{1}(p,\sigma)=h'(\sigma)p$ for some map 
$h':S\to H$, and $\pi_{L}\theta_{2}$ equals the projection $P\times S\to S$. Indeed, this is the case
iff $(\theta_{1},\theta_{2}):P\times S\to P\times L$ induces a section $P\times_{H}S\to P\times_{H}L$.

Since $(\pi_{L}:L\to S,\nabla^{L},h)$ is a prequantization for the action of $\mathfrak{h}$ on $(S,\omega^{S})$, 
the curvature of $\nabla^{L}$ equals $\omega^{S}$ and the representation of 
$\mathfrak{h}$ on $\Gamma^{\infty}(L)$ is given by Kostant's formula
\[X\mapsto \nabla^{L}_{\beta(X)}-2\pi\,i\la\mu,X\ra.\]

We identify $T(P\times_{H}S)\simeq (T P\times_{H} T S)/\mathord\sim$, cf.\ Example \ref{princ}.
For each $H$-equivariant vector field $v$ on $S$ and each $H$-equivariant vector field $w$ on $P$, let $[v,w]$ denote the 
vector field induced on $P\times_{H}S$.
For each $\theta=(\theta_{1},\theta_{2})\in\Gamma^{\infty}(P\times_{H}L)$ and $[w,v]\in\XX^{\infty}(P\times_{H}S)$, define
\[\nabla_{[w,v]}\theta:=\bigl(\theta_{1},\nabla^{L}_{v-\beta(\tau(w))}\theta_{2}-\tau_{h}(w)\cdot\theta_{2}\bigr),\]
where $(\nabla^{L}\theta_{2})(p,\sigma):=\nabla^{L}(\theta_{2}(p,\cdot))(\sigma)$ and $\tau_{h}(w)\in\mathcal{H}$ is the horizontal 
projection $w-\alpha(\tau(w))$ of $w$.
Suppose that $H$ is connected (we need this for equivariance of the connection $\nabla^{L}$, cf.\ Corollary \ref{eqcon}).
\begin{lemma}\label{concur}
$\nabla$ is a connection on $P\times_{H}L$ with curvature $\tilde{\omega}$.
\end{lemma}
\begin{proof}
First we check that $\nabla$ is well-defined. Indeed, $\nabla_{[\alpha(X),\beta(X)]}=0$, since
$\tau(\alpha(X))=X$ and $\tau_{h}(\alpha(X))=0$ and
\begin{align*}
h\cdot\nabla_{[w,v]}\theta&=\bigl(h\cdot\theta_{1},h\nabla^{L}_{v-\beta(\tau(w))}\theta_{2}- h\cdot(\tau_{h}(w))\cdot\theta_{2}\bigr)\\
&=\bigl(\theta_{1},\nabla^{L}_{h\cdot(v-\beta(\tau(w)))}h\cdot\theta_{2}-\tau_{h}(w)\cdot\theta_{2}\bigr)\\
&=\bigl(\theta_{1},\nabla^{L}_{v-\beta(\tau(w))}\cdot\theta_{2}-\tau_{h}(w)\cdot\theta_{2}\bigr),
\end{align*}
by $H$-equivariance of $\theta,\nabla,w, v,\beta,\alpha$ and $\tau$.

It is easy to check that $\nabla$ is a connection. For example, for an $H$-invariant function $f\in C^{\infty}(P\times S)^{H}$
one computes
\begin{align*}
\nabla_{[w,v]}f \theta&=\bigl(\theta_{1},\nabla^{L}_{v-\beta(\tau(w))}f\theta_{2})+\tau_{h}(w)f\theta_{2}\bigr)\\
&=\Bigl(\theta_{1},\bigl(f\nabla_{[w,v]}+\bigl(v-\beta(\tau(w))\bigr)\cdot f+\bigl(w-\alpha(\tau(w))\bigr)\cdot f\bigr)\theta_{2}\Bigr)\\
&=(f\nabla_{[w,v]}\theta+(w+v)f)\theta,\\
\end{align*}
since $(\alpha(\tau(w))+\beta(\tau(w)))\cdot f=0$ by invariance of $f$.

Now we shall compute the curvature of $\nabla$. Note the two ways in which the brackets $[,]$ are used, namely as a
commutator bracket and as a way to denote equivalence classes.
Let $[w_{1},v_{1}],[w_{2},v_{2}]$ be vector fields on $P\times_{H}S$ and $\theta=(\theta_{1},\theta_{2})\in\Gamma^{\infty}(P\times_{H}L)$.
We compute
\begin{align}\label{eqn1}
\bigl[\nabla_{[w_{1},v_{1}]},\nabla_{[w_{2},v_{2}]}\bigr]\theta
&=\Bigl(\theta_{1},
\bigl([\nabla^{L}_{v_{1}-\beta(\tau(w_{1}))},\nabla^{L}_{v_{2}-\beta(\tau(w_{2}))}]
+[\tau_{h}(w_{1}),\nabla^{L}_{v_{2}-\beta(\tau(w_{2}))}]\nonumber\\
&\quad+[\nabla^{L}_{v_{1}-\beta(\tau(w_{1}))},\tau_{h}(w_{2})]+[\tau_{h}(w_{1}),\tau_{h}(w_{2})]\bigr)\,\theta_{2}\Bigr)
\end{align}
Note that
\begin{align*}
[\tau_{h}(w_{1}),\nabla^{L}_{v_{2}-\beta(\tau(w_{2}))}]&=\nabla^{L}_{\tau_{h}(w_{1})\cdot\beta(\tau(w_{2}))}\\
&=\nabla^{L}_{\tau(w_{1})\cdot\beta(\tau(w_{2}))},
\end{align*}
and in the same way we obtain
\[[\nabla^{L}_{v_{1}-\beta(\tau(w_{1}))},\tau_{h}(w_{2})]=-\nabla^{L}_{\tau(w_{2})\cdot\beta(\tau(w_{1}))}.\]
We shall also need that
\begin{align*}
\lbrack\alpha(\tau(w_{1})),w_{2}]&=-\alpha(w_{2}\cdot\tau(w_{1}));\\
\lbrack w_{1},\alpha(\tau(w_{2}))]&=\alpha(w_{1}\cdot\tau(w_{2}));\\
\lbrack v_{1},\beta(\tau(w_{2}))]&=0;\\
\lbrack\beta(\tau(w_{1})),v_{2}]&=0,
\end{align*}
as follows from $H$-equivariance of $w_{1},w_{2},v_{1}$ and $v_{2}$.

On the other hand,
\begin{align}\label{eqnster}
\nabla_{\left[[w_{1},v_{1}],[w_{2},v_{2}]\right]_{T P\times T S}}\theta
&=\nabla_{\left[[w_{1},w_{2}]_{T P},[v_{1},v_{2}]_{T S}\right]}\theta\nonumber\\
&=(\theta_{1},\nabla^{L}_{[v_{1},v_{2}]_{T S}-\beta(\tau([w_{1},w_{2}]_{T P}))}\theta_{2}-\tau_{h}([w_{1},w_{2}]_{T P})\theta_{2}).
\end{align}
Using the well-known formula
\begin{align*}
F(w_{1},w_{2})&=d\tau(w_{1},w_{2})+[\tau,\tau]_{\mathfrak{h}}(w_{1},w_{2})\\
&=\tau([w_{1},w_{2}])-w_{1}\cdot(\tau(w_{2}))+w_{2}\cdot(\tau(w_{1}))+[\tau(w_{1},\tau(w_{2})]_{\mathfrak{h}},
\end{align*}
we continue the calculation of (\ref{eqnster})
\begin{align}\label{eqn2}
\nabla_{\left[[w_{1},v_{1}],[w_{2},v_{2}]\right]_{T P\times T S}}\theta&=\Bigl(\theta_{1},\bigl(\nabla^{L}_{[v_{1},v_{2}]_{T S}-\beta\left(F(w_{1},w_{2})+w_{1}\cdot\tau(w_{2})-w_{2}\cdot\tau(w_{1})
-[\tau(w_{1},\tau(w_{2})]_{\mathfrak{h}}\right)}\nonumber\\
&\quad+[w_{1},w_{2}]-\alpha(F(w_{1},w_{2})+w_{1}\cdot\tau(w_{2})-w_{2}\cdot\tau(w_{1})\nonumber\\
&\quad-[\tau(w_{1},\tau(w_{2})]_{\mathfrak{h}})\bigr)\,\theta_{2}\Bigr).
\end{align}
Note that $\theta_{2}$ is equivariant, hence for any $X\in\mathfrak{h}$ one has
\[X\cdot\theta_{2}=(\nabla^{L}_{\beta(X)}-2\pi i\la\mu,X\ra+\alpha(X))\theta_{2}=0\]
In particular, this is true for $X=F(w_{1},w_{2})$.

Subtracting the identity (\ref{eqn2}) from (\ref{eqn1}) one obtains the curvature, using all the given equalities, namely
\begin{align*}
K([w_{1},v_{1}],[w_{2},v_{2}])&=\frac{1}{2\pi i}\bigl([\nabla_{[w_{1},v_{1}]},\nabla_{[w_{2},v_{2}]}]-
\nabla_{\left[[w_{1},v_{1}],[w_{2},v_{2}]\right]_{T P\times T S}}\bigr)\\
&=\frac{1}{2\pi i}\Bigl([\nabla^{L}_{v_{1}-\beta(\tau(w_{1}))},\nabla^{L}_{v_{2}-\beta(\tau(w_{2}))}]\\
&\quad-\nabla^{L}_{[v_{1}-\beta(\tau(w_{1})),v_{2}-\beta(\tau(w_{2}))]_{T S}}
-2\pi i\la\mu,F(w_{1},w_{2})\ra\Bigr)\\
&=\omega\bigl(v_{1}-\beta(\tau(w_{1})),v_{2}-\beta(\tau(w_{2}))\bigr)-\la\mu,F(w_{1},w_{2})\ra\\
&=\tilde{\omega}([w_{1},v_{1}],[w_{2},v_{2}]).
\end{align*}
\end{proof}

A Hermitian metric $h'$ on $P\times_{H}L$ is given by 
\[h'([p,z_{1}],[p,z_{2}]):=h(z_{1},z_{2}).\]
This is well-defined since the representation of $H$ on $L$ is unitary.
\begin{lemma}
The connection $\nabla$ on $P\times_{H}L\to S'$ is Hermitian with respect to Hermitian metric $h'$.
\end{lemma}
\begin{proof}
This follows by computation: (in the notation introduced previously and $\theta=(\theta_{1},\theta_{2})$,
$\theta'=(\theta_{1}',\theta_{2}')$)
\begin{align*}
h'(\nabla_{[w,v]}\theta,\theta')+h'(\theta,\nabla_{[w,v]}\theta')&=h(\nabla^{L}_{v-\beta(\tau(w))}\theta_{2},\theta_{2}')
+h(\tau_{h}(w)\cdot\theta_{2},\theta_{2}')\\
&\quad + h(\theta_{2},\nabla^{L}_{v-\beta(\tau(w))}\theta_{2}')
+h(\theta_{2},\tau_{h}(w)\cdot\theta_{2}')\\
&=(v-\beta(\tau(w)))\cdot h(\theta_{2},\theta_{2}')+\tau_{h}(w)\cdot h(\theta_{2},\theta_{2}')\\
&=[w,v]\cdot h'(\theta,\theta'),
\end{align*}
where in the third line we used the fact that $\nabla^{L}$ is Hermitian and in the last that $\theta_{1}$, $\theta_{2}$,
$\theta_{1}'$ and $\theta_{2}'$ are $H$-equivariant.
\end{proof}
\end{example}
\begin{corollary}
The triple $(P\times_{H}L\to P\times_{H}S,\nabla,h')$ is a prequantization for the action of $P\times_{H}P$ on 
$(P\times_{H}S,\tilde{\omega})$, with prequantization representation 
\[\AA(P\times_{H}P)\simeq T P/H\to \mathcal{D}(P\times_{H}L)\] 
given by
\begin{align*}
w&\mapsto\nabla_{\gamma(w)}-2\pi i \la\tilde{\mu},w\ra\\
&\quad=\nabla_{[w,0]}-2\pi i\la\mu,\tau(w)\ra.
\end{align*}
\end{corollary}

Note that in Lemma \ref{concur} we have used the $H$-equivariance of the connection $\nabla^{L}$. We shall now prove a more 
general result for source-connected Lie groupoids.
\begin{lemma}\label{coneq}
Consider the situation of Theorem \ref{preqrep}. For any prequantization representation of a $\AA$ on a line bundle $L\to S$, the given connection 
$\nabla$ on $L$ is $\AA$-equivariant.
\end{lemma}
\begin{proof}
One computes for any $v\in\XX_{\tilde{J}}^{\infty}(S)$
\begin{align*}
[\pi(X),\nabla_{v}]&=[\nabla_{\alpha(X)},\nabla_{v}]+2\pi i\la\tilde{\mu},J^{*}X\ra\\
&=2\pi i \omega(\alpha(X),v)+\nabla_{[\alpha(X),v]}+2\pi i v\cdot\la\mu,J^{*}X\ra\\
&=\nabla_{[\alpha(X),v]},
\end{align*}
which means exactly that $\nabla$ is $\AA$-equivariant.
\end{proof}
The corollary that we tacitly used in the proof of Lemma \ref{concur} (in the particular case that $G$ is a Lie group) is
\begin{corollary}\label{eqcon}
If $G$ is a source-connected integrating Lie groupoid of $\AA$, then $\nabla$ is equivariant, in the sense 
that for $v\in T^{J}_{\sigma}S$ and each $g\in G_{J(\sigma)}$
\[\nabla_{g v}=g\nabla_{v}g^{-1}.\]
\end{corollary}
\begin{proof}
Choose a connection $\nabla'$ on $\AA$. Then there exists an exponential map $\exp_{\nabla'}:\AA\to G$ (cf.\ \cite{La})
Differentiating the expression 
\[\nabla_{\exp_{\nabla'}(\tau X) v}=\exp_{\nabla'}(\tau X)\nabla_{v}\exp_{\nabla'}(-\tau X).\]
at $\tau=0$ gives the equality in the proof above. 
\end{proof}

\subsection{Integrating a prequantization representation of a Lie algebroid.}\label{integ}
In this section we discuss the integrability of Lie algebroid representations. In particular, we consider the examples from the previous section.

Not every Lie algebroid integrates to a Lie groupoid. Precise conditions for the existence of an integrating Lie groupoid
for a given Lie algebroid are given in \cite{CF}. Suppose $\AA$ is a Lie algebroid and $\alpha':\AA\ltimes J\to\mathcal{D}_{U}(L)$ 
a Hermitian representation (e.g., obtained by prequantization). One would like to integrate such a representation
to a representation of a Lie groupoid which has associated Lie algebroid $\AA$. 
\begin{definition}
A \textbf{representation of a Lie groupoid $G\rra M$} on a smooth complex vector bundle $E\to M$ is a smooth 
action of $G$ on $\pi:E\to M$
\[\pi:G\fiber{s}{\pi}E\to E\]
that is linear, i.e.
\[\pi(g,\lambda\cdot e)=\lambda\cdot\pi(g,e)\]
and
\[\pi(g,e+f):=\pi(g,e)+\pi(g,f)\]
for all $g\in G$, $\lambda\in\CC$ and $e,f\in E_{s(g)}$.
\end{definition}
\begin{remark}
This notion generalizes the notion of $H$-equivariant vector bundle for a Lie group $H$.
\end{remark}
The representation $\pi$ is \textbf{unitary} with respect to a Hermitian metric $h$ on $E$ if it preserves $h$, i.e.
\[h(\pi(g,e),\pi(g,f))=h(e,f),\]
for all $g\in G$ and $e,f\in E_{s(g)}$.
A unitary representation $\pi$ can equivalently be given by a morphism of groupoids $G\to U(E)$, where $U(E)$ is the Lie groupoid of linear 
unitary maps $E_{m}\to E_{n}$ for all $m,n\in M$ (it has a natural smooth structure, cf.\ \cite{Mac}).
\begin{remark}
More generally one could consider continuous unitary representations of groupoids on continuous fields of Hilbert spaces.
These are studied in \cite{Bos}.
\end{remark}
Suppose $G$ acts on a map $J:N\to M$. Suppose $E\to N$ a complex smooth vector bundle endowed with Hermitian
structure $h$.
\begin{definition}
A \textbf{(unitary)representation of $G$} on a smooth complex vector bundle $E\to N$ is  
a (unitary) representation of the action groupoid $G\ltimes J$ on $\pi:E\to N$.
\end{definition}

Suppose $\AA$ is integrable and $J:S\to M$ is proper, then by Proposition 3.5 and Proposition 5.3 in \cite{MM2} 
the representation $\pi:\AA\ltimes J\to \mathcal{D}(E)$ of the Lie algebroid $\AA$ on the vector bundle $E\to S$
integrates to a unitary representation $G\ltimes J\to U(L)$ of the source-simply connected integrating Lie 
groupoid $G$ of $\AA$ on $E\to S$. The condition that $J$ is proper will also arise in the next section about
the quantization procedure. Note that one can prove that a proper smooth family of manifolds is a fiber bundle.

\begin{example}
A flat connection $\nabla-2\pi i\mu$ on a line bundle $L\to M$ is a prequantization representation of $T M$ as a Lie algebroid
acting Hamiltonianly on $(M,\omega)$ as in Example \ref{flatcon}. It integrates to a representation of a source-simply connected
Lie groupoid integrating $T M$, for example the fundamental groupoid $\pi_{1}(M)$ of $M$.
The representation is exactly the parallel transport associated to the connection $\nabla-2\pi i\mu$.

The only prequantization that lifts to a representation of the pair groupoid $M\times M$ (which also integrates $T M$)
is the representation $d-2\pi i\mu$ on the trivial line bundle $M\times\CC\to M$.
\end{example}

\begin{example}
Recall the situation of Example \ref{princ4}.
There exists a canonical unitary representation $\bar{\pi}$ of $P\times_{H}P$ on $P\times_{H}L$, defined by
\[\bar{\pi}([p,q])[q,z]:=[p,z]\]
for suitable representatives.
\begin{proposition}\label{princ5}
The representation $\bar{\pi}$ integrates the prequantization representation given by 
\[w\mapsto \nabla_{\gamma(w)}-2\pi\,i\la\tilde{\mu},w\ra\]
(cf.\ Theorem \ref{preqrep}), where $\nabla$ is the connection we have defined  in Lemma \ref{concur} and\linebreak
\mbox{$\tilde{\mu}:P\times_{H}S\to T^{*}P/H$} is the momentum map for the Hamiltonian action of $P\times_{H}P\rra M$ on 
$P\times_{H}S\to M$ (cf.\ Proposition \ref{princ6}).
\end{proposition}
Note therefore, that $\bar{\pi}$ does not depend on the chosen connection $\tau$.
\begin{proof}
Suppose that $w\in\XX^{\infty}(P)$ is $H$-equivariant. It represents a smooth section in 
$\Gamma^{\infty}(\AA(P\times_{H}P))\simeq\Gamma^{\infty}(T P/ H)$.
Note that the infinitesimal representation $\bar{\pi}':\AA(P\times_{H}P)\to\mathcal{D}(L)$ associated to the representation
$\bar{\pi}$ is given by
\[\bar{\pi}'(w)\theta=(\theta_{1},w\theta_{2}).\]
On the other hand,
\begin{align*}
\pi(w)\theta&=\nabla_{\gamma(w)}\theta-2\pi i\la\tilde{\mu},w\ra\theta\\
&=\nabla_{[w,0]}\theta-2\pi i\la\tilde{\mu},w\ra\theta\\
&=(\theta_{1},\nabla^{L}_{-\beta(\tau(w))}\theta_{2}+(w-\alpha(\tau(w)))\theta_{2}-2\pi i\la\mu,\tau(w)\ra\theta_{2})\\
&=(\theta_{1},w\cdot\theta_{2}),
\end{align*}
from which the lemma follows.
\end{proof}
\end{example}

\section{Quantization and symplectic reduction}
\subsection{Quantization through K\"{a}hler polarization.}
In this section we introduce K\"{a}hler quantization for Hamiltonian Lie algebroid actions. Next, we discuss the examples
which we have been considering throughout the paper.

Suppose a regular Lie algebroid $(p:\AA\to M,\rho)$ acts in a Hamiltonian fashion on a smooth surjective submersion
$J:S\to M$, with a $J$-presymplectic 2-form 
$\tilde{\omega}\in\Omega_{\tilde{J}}^{2}(S)$, where $\tilde{J}:S\to M/\rho$. 
Denote the action by $\alpha:\Gamma^{\infty}(\AA)\to\mathfrak{X}^{\infty}_{\tilde{J}}(S)$.
Let $\tilde{\mu}:S\to J^{*}\AA^{*}$ be a momentum map.

In this section we shall make the additional assumption that $J:S\to M$ is a smooth bundle of compact connected K\"{a}hler manifolds.
Denote the almost complex structure by
\[j:T^{J}S\to T^{J}S.\] 
The following conditions are satisfied: $\omega(j\cdot,\cdot)>0$ and $\omega(j\cdot,j\cdot)=\omega$,  where $\omega:=\tilde{\omega}|_{T^{J}S}$.

Our final assumption is that the almost complex structure $j$ is $\AA$-equivariant, in the sense that
\[[\alpha(X), j(v)]_{T^{\tilde{J}}S}=j[\alpha(X), v]_{T^{\tilde{J}}S},\]
for all $X\in\Gamma^{\infty}(\AA)$ and $v\in \mathfrak{X}_{J}^{\infty}(S)$.

Let $T^{J,\CC}S$ denote the complexification $T^{J}S\otimes\CC$ of $T^{J}S$. 
The complex extension of $j$ is denoted by
$j_{\CC}:T^{J,\CC}S\to T^{J,\CC}S$. 
\begin{definition}
The \textbf{K\"ahler polarization} $\mathcal{P}(S,j)$ of $(J:S\to M,\omega)$ is defined by
\[\mathcal{P}(S,j):=\{v\in T^{J,\CC}S\mid j_{\CC}(v)=-i v\in T^{J,\CC}S\}.\]
\end{definition}
Smooth sections of $\mathcal{P}(S,j)$ are called \textbf{polarized}.

Assume $(L\to S,\nabla^{L},h)$ is a prequantization line bundle for the action of $\AA$. Denote the associated representation 
(see \S\ref{preq}) by $\pi:\AA\ltimes J\to\mathcal{D}(L)$.

\begin{definition}
The \textbf{geometric quantization} of a prequantization line bundle \linebreak \mbox{$(L\to S,\nabla,h)$} for the Hamiltonian action $\alpha$ of $\AA$
on $(J:S\to M,\tilde{\omega})$ is given by
\[\Delta_{Q}:=\{\sigma\in\Gamma^{\infty}(L)\mid\nabla_{v}\sigma=0\mbox{ for all }v\in\Gamma^{\infty}(\mathcal{P}(S,j))\}.\]
\end{definition}
We call the sections of $L\to S$ in $\Delta_{Q}$ \textbf{holomorphic}.
The space $\Delta_{Q}$ has the structure of a Hilbert $C_{0}(M)$-module (cf.\ \cite{Bos}), with
$C_{0}(M)$-valued inner product
\[\la\sigma,\sigma'\ra(m):=\int_{S_{m}}h(\sigma,\sigma')\Omega_{m},\]
where $\{\Omega_{m}\}_{m\in M}$ is a smooth family of densities on $J:S\to M$ defined by
\[\Omega_{m}:=\omega_{m}^{d_{m}}/(d_{m}!),\]
($d_{m}:=\dim(S_{m})/2$).
Hence, it corresponds to a continuous field of Hilbert spaces over $M$. 
We shall assume that it actually is a smooth vector bundle, which we denote by $Q\to M$.
\begin{theorem}\label{qu}
The geometric quantization vector bundle $Q\to M$ carries a Hermitian representation of $\AA$.
\end{theorem}
\begin{proof}
The fiber of the vector bundle $Q\to M$ at $m$ is given by 
\[\{\sigma\in\Gamma^{\infty}(S_{m},L_{m})\mid\nabla_{v}\sigma=0\mbox{ for all }v\in\Gamma^{\infty}(\mathcal{P}(S,j)|_{S_{m}})\},\]
which is finite-dimensional, since $S_{m}$ is compact.

We check that the representation $\pi$ of $\AA$ on $L\to M$ restricts to $\Delta_{Q}$. From this the theorem follows. Suppose $\nabla_{v}\sigma=0$
for all $v\in\Gamma^{\infty}(\mathcal{P}(S,j))$. Note that $\AA$-equivariance of $j$ implies that $[\alpha'(X),v]\in \mathcal{P}(S,j)$ whenever $v\in \Gamma^{\infty}(\mathcal{P}(S,j))$.
Indeed, suppose $v\in\Gamma^{\infty}(\mathcal{P}(S,j))$ and $X\in\Gamma^{\infty}(\AA)$ then
\begin{align*}
j_{\CC}[\alpha'(X),v]&=[\alpha'(X),j_{\CC}v]\\
&=[\alpha'(X),-i v]\\
&=-i [\alpha'(X),v].
\end{align*}
Hence, the $\Delta_{Q}(S,\tilde{\omega})$ is $\AA$-invariant: for $\sigma\in\Delta_{Q}(S,\tilde{\omega})$
\begin{align*}
\nabla_{v}(\pi(X)\sigma)&=\nabla_{v}(\nabla_{\alpha'(X)}\sigma-2\pi i\la\tilde{\mu},J^{*}X\ra\sigma)\\
&=\nabla_{\alpha'(X)}\nabla_{v}\sigma-\nabla_{[\alpha'(X),v]}\sigma-2\pi i\, \tilde{\omega}(\alpha'(X),v)\sigma\\
&\quad+2\pi i\la\tilde{\mu},J^{*}X\ra\nabla_{v}\sigma-2\pi i (v\cdot\la\tilde{\mu},J^{*}X\ra)\sigma\\
&=-2\pi i v\inner(d^{J}\la\tilde{\mu},J^{*}X\ra+\alpha'(X)\inner\tilde{\omega})\sigma\\
&=0,
\end{align*}
by the quantization condition (\ref{quancond}).

One easily sees that the representation is Hermitian since the prequantization representation is Hermitian and $\omega$ is invariant.
Indeed, for all $m\in M$,$\sigma,\sigma'\in\Gamma^{\infty}(L)$ and $X\in\Gamma^{\infty}(\AA)$
\begin{align*}
(\la\pi(X)\sigma,\sigma'\ra+\la\sigma,\pi(X)\sigma'\ra)(m)&=\int_{S_{m}}h(\pi(X)\sigma,\sigma')+h(\sigma,\pi(X)\sigma')\Omega_{m}\\
&=\int_{S_{m}}\alpha(X)\cdot h(\sigma,\sigma')\Omega_{m}\\
&=\bigl((T J\circ\alpha(X))\cdot(m'\mapsto \int_{S_{m'}}h(\sigma,\sigma')\Omega_{m'})\bigr)(m),\\
\end{align*}
which finishes the proof.
\end{proof}

\begin{remark}
Of course, one might wonder whether one can integrate this representation to a representation of an integrating Lie 
groupoid. This issue is somewhat simpler than for prequantization line bundles (cf.\ section \ref{integ}), since
$Q\to M$ is a vector
bundle \textit{over} $M$. Hence, if $G$ is a source-connected integrating Lie groupoid for $\AA$, then every quantization representation  
of $\AA$ integrates to a representation of $G$.
\end{remark}

\begin{example}
In the case of the Hamiltonian action of an integrable distribution\linebreak \mbox{$T \mathcal{F}\subset T M\to M$} on $(J:M\to M,\omega)$, with momentum map $\mu:M\to \mathcal{F}^{*}$, a prequantization representation on a line bundle $L\to M$ with metric $h$ and Hermitian connection $\nabla$
is given by $\nabla-2\pi i\mu$, where $d^{\mathcal{F}}\mu=-\omega$. Obviously, the quantization procedure is empty in this
situation, since the fibers of $J$ are points.
\end{example}

\begin{example}
If $\mathfrak{g}\to M$ is a smooth family of Lie algebras that acts in a Hamiltonian fashion on a bundle of symplectic manifolds $S\to M$ and there is a prequantization
$(L,\nabla,h)$, then K\"{a}hler quantization is family K\"{a}hler quantization.
\end{example}

\begin{example}\label{princ8}
Suppose $H$ is a Lie group that acts in a Hamiltonian fashion on a symplectic manifold $(S,\omega_{S})$ with momentum map $\mu$. Suppose 
$(L,\nabla^{L},h)$ is a prequantization of this action. Furthermore, suppose that $P$ is a principal $H$-bundle,
endowed with an $\mathfrak{h}$-valued connection 1-form $\tau$.
In Example \ref{princ3} it was shown that there exists a closed form $\tilde{\omega}$ on $S':=P\times_{H}S$, such that
the action of the gauge groupoid $P\times_{H}P$ on $(S',\tilde{\omega})$ is Hamiltonian. 
In Example \ref{princ4} it was shown that there exists a prequantization $(P\times_{H}L\to S',\nabla,h')$ of 
this action.

Suppose that $S$ is a compact K\"{a}hler manifold with $H$-equivariant almost complex structure $j:T S\to T S$.
\begin{lemma}\label{complstr}
The almost complex structure $j$ induces a $P\times_{H}P$-equivariant family of almost complex structures
\[j':T^{J}S'\to T^{J}S'\]
on $J:S'\to M$. 
\end{lemma}
\begin{proof}
We shall use the isomorphism $T^{J}(P\times_{H}S)\simeq P\times_{H} T S$ from Proposition \ref{princ2}. Define the almost complex 
structure as a map $j':P\times_{H} T S\to P\times_{H} T S$ by
\[j'([p,v])=[p,j(v)].\]
This is obviously a complex structure:
\[j'(j'([p,v]))=[p, j(j(v))]=[p,-v]=-[p,v],\]
which is $P\times_{H}P$-equivariant by the computation
\[j'([p,q]\cdot[q,v])=j'([p,v])=[p,j(v)]=[p,q]\cdot[q,j(v)]=[p,q]\cdot j'([q,v]).\]
\end{proof}
From the lemma we conclude that  $J:S'\to M$ is a bundle of K\"ahler manifolds.
So there exists a K\"{a}hler quantization vector bundle $Q'\to M$ of the prequantization 
\[(P\times_{H}L\to S',\nabla,h')\]
of the action of $P\times_{H}P$ on $(S'\to M,\tilde{\omega})$.

Let $Q$ denote the representation space obtained by quantization of the action of $H$ on $S$ with 
prequantization $(L,\nabla^{L},h)$.
The associated vector bundle $P\times_{H}Q\to M$ carries a canonical representation of $P\times_{H}P$.
\begin{proposition}\label{quaniso}
The vector bundle $P\times_{H}Q\to M$ is $P\times_{H}P$-equivariantly isomorphic to the quantization bundle 
$Q'\to M$.
\end{proposition}
\begin{proof}
Note that by definition there is a bijection between sections of $Q'\to M$ and holomorphic sections of 
$P\times_{H}L\to S'$. 
\begin{lemma}
There exists a canonical fiberwise linear bijection
\[P\times_{H}\Gamma^{\infty}(L)\to \Gamma^{\infty}(P\times_{H}L),\]
where $\Gamma^{\infty}(P\times_{H}L)$ is thought of as a vector bundle with fiber at $m$ given by
\[\Gamma^{\infty}(P\times_{H}L|_{P_{m}\times_{H}S}).\]
\end{lemma}
\begin{proof}
A bundle morphism $\Psi:P\times_{H}\Gamma^{\infty}(L)\to \Gamma^{\infty}(P\times_{H}L)$ is defined by
\[(m, [p,\eta])\mapsto \Bigl(m,\bigl([p',\sigma]\mapsto\bigl[p,\eta\bigl((p (p')^{-1})\cdot\sigma\bigr)\bigr]\bigr)\Bigr),\]
where $m\in M$, $p\in\pi^{-1}(m)$, $\eta\in\Gamma^{\infty}(L)$ and $p' p^{-1}$ is the unique element in $H$ such that 
$(p' p^{-1})\cdot p=p'$.
This is well-defined. Indeed, for fixed $p\in P$ one has
\begin{align*}[h\, p',h\, \sigma]&\mapsto \bigl[p,\eta\bigl(((p(h\,p')^{-1})\cdot h \sigma\bigr)\bigr]\\
&=\bigl[p,\eta\bigl(((p (p')^{-1})\cdot \sigma\bigr)\bigr]
\end{align*}
and
\begin{align*}[h\,p,h\cdot\eta]&\mapsto\Bigl([p',\sigma]\mapsto \bigl[h\,p,h\,\eta(h^{-1}((h\, p)(p')^{-1})\cdot\sigma)\bigr]\Bigr)\\
&=([p',\sigma]\mapsto [p,\eta((p (p')^{-1})\cdot\sigma)]).
\end{align*}
The map $\Psi$ is obviously linear. A two-sided inverse is as follows. Suppose $(\theta_{1},\theta_{2})$ is a
section of $P_{m}\times_{H}L\to P_{m}\times_{H}S$ (cf.\ \ref{princ4}). Define a map $\Phi:\Gamma^{\infty}(P\times_{H}L)\to P\times_{H}\Gamma^{\infty}(L)$ by
\[[m,(\theta_{1},\theta_{2})]\mapsto\Bigl(\sigma\mapsto\bigl(m,\bigl[p,(p(\theta_{1}(p,\sigma))^{-1})\cdot\theta_{2}(p,\sigma)\bigr]\bigr)\Bigr),\]
for a chosen $p\in P$. Straightforward calculations using the equivariance of $\theta_{1}$ and $\theta_{2}$
show that this is independent of the choice of $p\in P_{m}$ and that $\Psi\circ\Phi=1$ and $\Phi\circ\Psi=1$.
\end{proof}
On sections one obtains, for a smooth map $(m,\sigma)\mapsto\eta_{m}(\sigma)$ from $M\times S$ to $Q$ and a 
smooth section $\xi\in\Gamma^{\infty}(P)$, a
smooth section $\Psi(\xi,\eta)\in\Gamma^{\infty}(P\times_{H}L)$ given by
\[[p,\sigma]\mapsto[\xi(\pi(p)),\eta_{\pi(p)}((\xi(\pi(p))p^{-1})\cdot\sigma)].\]

One easily checks that the sections $\eta_{m}\in\Gamma^{\infty}(L)$ are holomorphic for all $m\in M$ iff 
the image $\Psi([\xi,\eta])$ is holomorphic. 
Indeed, for 
$[\xi',v']\in \mathcal{P}(S',j')$, with $\xi':P\times S\to P$ and $v':P\times S\to T^{\CC}S$ 
such that $v'(p)$ is polarized for each $p\in P$, we have
\begin{align*}
\nabla_{[\xi',v']}\Psi(\xi,\eta)(p,\sigma)&=[\xi(p,\sigma),\nabla_{(\xi\xi'^{-1})v'(p,\cdot)}(\eta\circ
\beta(\xi(\pi(p))p^{-1})(\sigma)]\\
&=[\xi(\pi(p)),0]\\
&=0,
\end{align*}
by equivariance of $\nabla^{L}$. The reverse statement is proven by the same formula.
\end{proof}
\end{example}

\begin{remark}\label{KKstuff}
We shall sketch a more general view of geometric quantization based on \cite{La3} and \cite{GGK} and references in these papers.
If the line bundle $L\to S$ is positive enough, then the quantization $Q\to M$ vector bundle equals the index of
a continuous $G$-equivariant family of Dolbeault-Dirac operators 
\[\{{\bar\partial}_{L_{m}}+{\bar\partial}^{*}_{L_{m}}:\Omega^{0,\mathrm{even}}(S_{m};L_{m})\to\Omega^{0,\mathrm{odd}}(S_{m};L_{m})\}_{m\in M}\] 
constructed from the connection $\nabla^{L}$ on the line bundle $L\to S$ and the family of almost complex
structures on $S$.

For now, suppose that $G$ is locally compact, $\sigma$-compact, endowed with a Haar system and that the action of $G$ on $S$ is proper.
Then this family of Dolbeault-Dirac operators gives rise to a cycle in Kasparov's $G$-equivariant bivariant K-theory
\[[\{{\bar\partial}_{L_{m}}+{\bar\partial}^{*}_{L_{m}}\}_{m\in M}]\in KK_{G}(C_{0}(S),C_{0}(M)).\]
The Baum-Connes analytical assembly map for groupoids (cf.\ e.g., \cite{Con}, \cite{Tu} or \cite{Pat})
\[\mu:KK_{G}(C_{0}(S),C_{0}(M))\to KK(\CC,C^{*}_{r}(G))\simeq K_{0}(C^{*}_{r}(G))\]
maps the class $[D]:=[\{{\bar\partial}_{L_{m}}+{\bar\partial}^{*}_{L_{m}}\}_{m\in M}]$ to a class in the K-theory of the
$C^{*}$-algebra $C^{*}_{r}(G)$ of the groupoid $G$.

A different way to look at geometric quantization is to define $\mu([D])$ to be the geometric quantization of the Hamiltonian action of $G\rra M$
on $(J:S\to M,\tilde{\omega})$.
Under certain conditions, including $G$ being a proper groupoid, $K_{0}(C^{*}_{r}(G))$ is isomorphic to the
representation ring of $G$ (cf.\ \cite{TuXu}). Hence, geometric quantization in this sense will then still yield
(a formal difference of) representations of $G$.
This approach gives new possibilities to generalize geometric quantization. Instead of requiring $J:S\to M$ to be 
a bundle of compact K\"{a}hler
manifolds, one requires $J:S\to M$ to be endowed with a $G$-equivariant family of $\spin^{c}$-structures. 
One then proceeds by defining the geometric quantization
as the image under the analytical assembly map of KK-cycle defined by the associated family of $\spin^{c}$-Dirac operators coupled to the 
prequantization line bundle. This generalizes the notion of family quantization (cf. \cite{Zha}).
\end{remark}

\subsection{Symplectic reduction.}
In this section we discuss a generalization of symplectic reduction to our setting.
We reduce in stages, first internal symplectic reduction, then the entire symplectic reduction.

Suppose $G$ is a source-connected regular Lie groupoid. Suppose $\alpha$ is an {\em internally} strongly Hamiltonian left action of $G$ on a smooth 
bundle of connected symplectic manifolds \mbox{$(J:S\to M,\omega\in\Omega_{J}^{2}(S))$} with internal momentum map $\mu:S\to\AA^{*}(I_{G})$
(see Definition \ref{locstr}). 

Denote the image of the zero section $0:M\to\AA^{*}(I_{G})$ by $0_{M}$. Suppose $0_{M}\subset\im(\mu)$ and $\mu$ and $0$ are transversal, i.e.
$T 0(T M)$ and $T\mu(T S)$ are transversal in $T \AA^{*}(I_{G})$. Then $\mu^{-1}(0_{M})$ is a manifold. Suppose, furthermore, that $G_{m}^{m}$ acts freely and properly on $\mu^{-1}(0(m))$ for each $m\in M$.
Then for each $m\in M$ the quotient manifold 
\[S_{m}^{(0)}:=G_{m}^{m}\backslash\mu^{-1}(0)\] 
is a smooth manifold with a symplectic 2-form $\omega_{m}^{0}\in\Omega^{2}(S_{m})$ uniquely determined by the equation
\[p_{m}^{*}\omega^{0}_{m}=i^{*}_{m}\omega|_{S_{m}},\]
(cf.\ \cite{MaWe}).
\begin{lemma}
The map
\[\bigcup_{m\in M}S_{m}^{(0)}=I_{G}\backslash\mu^{-1}(0_{M})\to M\]
is a smooth family of symplectic manifolds.
\end{lemma}
\begin{definition}
The smooth family of symplectic manifolds 
\[I_{G}\backslash\mu^{-1}(0_{M})\to M\] 
is called the \textbf{internal Marsden-Weinstein quotient} of
the internally Hamiltonian groupoid action.
\end{definition}

\begin{example}
In the case of a smooth family of Lie groups acting on a smooth family of symplectic manifolds the internal Marsden-Weinstein
quotient is simply the family bundle of Marsden-Weinstein quotients of the actions of the groups on the fibers (cf. \cite{Zha}).
\end{example}
\begin{example}\label{MWprinc}
Suppose $H$ is a Lie group acting on a symplectic manifold $(S,\omega^{S})$ in a Hamiltonian fashion, with momentum
map $\mu$. Suppose $P$ is a principal $H$-bundle. 
Then $G:=P\times_{H}P\rra M$ acts on $(P\times_{H}S,\omega)$ in a internally Hamiltonian fashion (cf.\ Example \ref{princ2}), with 
momentum map $\bar{\mu}[p,s]:=[p,\mu(s)]\in P\times_{H}\mathfrak{h}^{*}$.
The internal Marsden-Weinstein quotient $I_{G}\backslash\bar{\mu}^{-1}(0_{M})$ is symplectomorphic to $M\times\mu^{-1}(0)/H$ as
a smooth bundle of symplectic manifolds, using the map 
\[I_{G}\backslash\bar{\mu}^{-1}(0_{M})\To M\times\mu^{-1}(0)/H,[p,\sigma]\mapsto (\pi(p),[\mu(\sigma)]).\]
\end{example}

We now turn our attention to the entire quotient.
\begin{lemma}
The map
\[G\backslash\mu^{-1}(0_{M})\to G\backslash M\]
can be given the structure of a continuous family of symplectic manifolds.
\end{lemma}
\begin{remark}
Note that the space $G\backslash M$ is in general non-Hausdorff, unless $G\rra M$ is proper. In the case it is Hausdorff, there is a smooth structure on $G\backslash\mu^{-1}(0_{M})$ and $G\backslash M$ if the action of $I_{G}\rra M$ on $J:S\to M$ is proper and free, as mentioned above,
$\mu$ and $0$ are transversal and the groupoid $G\rra M$ is effective.
\end{remark}
\begin{proof}
Since the momentum map $\mu:S\to\AA^{*}(I_{G})$ is equivariant and since for every $g\in G_{m}^{n}$ one has 
$\Ad^{*}(g)(0(m))=0(n)$, the smooth isomorphism $\alpha(g):S_{m}\to S_{n}$ restricts to a smooth isomorphism 
\[\alpha(g):\mu^{-1}(0(m))\to\mu^{-1}(0(n)).\]

This induces a well-defined action $\bar{\alpha}$ of $G$ on the internal Marsden-Weinstein quotient $I_{G}\backslash\mu^{-1}(0_{M})\to M$ given by
\[\bar{\alpha}(g)(G_{m}^{m}\sigma):=(G_{n}^{n}\alpha(g)\sigma).\]
Indeed, suppose $\sigma_{1}=\alpha(g')\sigma_{2}$, for $\sigma_{1},\sigma_{1}\in\mu^{-1}(0(m))$ and $g'\in G_{m}^{m}$. 
Then $\alpha(g)\sigma_{1}=\alpha(g\,g'\,g^{-1})\alpha(g)\sigma_{2}$.
Actually, $\bar{\alpha}(g)=\bar{\alpha}(h)$ for all $g,h\in G_{m}^{n}$, as one checks by a similar computation.

Of course the action $\bar{\alpha}$ induces an equivalence relation on $I_{G}\backslash\mu^{-1}(0_{M})$ and the quotient equals
$G\backslash\mu^{-1}(0_{M})$. For every $g\in G$, the map $\bar{\alpha}(g)$ is a symplectomorphism, since $\alpha(g)$ is a 
symplectomorphism. Hence there exists a canonical family of symplectic forms $\omega^{00}$ on $G\backslash\mu^{-1}(0_{M})\to G\backslash M$.
\end{proof}
\begin{definition}
The continuous bundle of symplectic manifolds $(G\backslash\mu^{-1}(0_{M})\to G\backslash M,\omega^{00})$ is called the \textbf{Marsden-Weinstein quotient} of the internally Hamiltonian action of $G\rra M$ on $(J:S\to M,\omega)$.
\end{definition}
\begin{example}
Consider a group $H$ that acts on a manifold $M$. Denote the action by $\alpha$. The action groupoid $H\ltimes M\rra M$ acts in a Hamiltonian
fashion on $(id:M\to M,0)$, with momentum map given by any $\mu:M\to(\mathfrak{h}\ltimes M)^{*}$ such that 
$d^{\mathfrak{h}\ltimes M}\mu=d(\alpha^{*}\mu)=0$. The Marsden-Weinstein quotient is defined iff $\mu=0$ and then, obviously, is given by $M/H\to M/H, 
0)$. Note that $M/H$ only is smooth if the action of $H$ on $M$ is proper and free.
\end{example}
\begin{example}
We continue Example \ref{MWprinc}. One easily sees that
the Marsden-Weinstein quotient $((P\times_{H} P)\backslash\tilde{\mu}^{-1}(0_{M})\to (P\times_{H} P)\backslash M=*,\tilde{\omega}^{00})$ is symplectomorphic to the
Marsden-Weinstein quotient $(H\backslash\mu^{-1}(0),(\omega^{S})^{0})$.
\end{example}

\subsection{Quantization commutes with reduction.}
Suppose $G\rra M$ is a regular Lie groupoid.
Suppose $\pi:G\to U(E)$ is a unitary representation of $G$ on a vector bundle $p:E\to M$ with hermitian structure $h$.
Define the continuous family of vector spaces $E^{I_{G}}\to M$ of $I_{G}$-fixed vectors by
\[E^{I_{G}}:=\{e\in E\mid \pi(g)e=e\mbox{ for all }g\in G_{p(e)}^{p(e)}\}.\]
We shall assume in this paper that it is actually a smooth vector bundle.
One easily checks that $E^{I_{G}}$ is closed under $G$. Indeed, for $g\in G_{n}^{n}$, $h\in G_{m}^{n}$ and $e\in E^{I_{G}}_{m}$ one has
\begin{align*}
\pi(g)\pi(h)e&=\pi(g h)e\\
&=\pi(h)\pi((h^{-1}g h))e\\
&=\pi(h)e.
\end{align*}
\begin{definition}
The smooth vector bundle $E^{I_{G}}\to M$ is called the \textbf{internal quantum reduction} of $\pi:G\to U(E)$.
\end{definition}
By similar reasoning as above, the restriction of $\pi$ to a map $G\to U(E^{I_{G}})$ is $I_{G}$-invariant, i.e.
\[\pi(g)e=\pi(g')e\]
for all $g,g'\in G_{m}^{n}$, $n,m\in M$ and $e\in E^{I_{G}}_{m}$. 

Suppose $G$ acts in Hamiltonian fashion on a smooth bundle of compact connected K\"{a}hler manifolds $J:S\to M$ endowed 
with a $J$-presymplectic form $\tilde{\omega}$, such that the complex structure $j$ and Hermitian metric $h$ are
$G$-equivariant.
Denote the momentum map by $\tilde{\mu}: S\to J^{*}\AA^{*}$.
Consider the internal Marsden-Weinstein quotient $(\mu^{-1}(0_{M})/I_{G},\omega^{0})$.
Suppose $(L,\nabla^{L},h)$ is a prequantization of the $G$-action. Suppose $L^{0}\to \mu^{-1}(0_{M})/I_{G}$ is a line bundle
such that
\[p^{*}L^{0}=L|_{\mu^{-1}(0_{M})}.\]
This is a strong condition, which is satisfied if the action of $I_{G}$ on $S$ is free.
The line bundle $L^{0}$ has an induced prequantization connection $\nabla^{0}$, since $\nabla^{L}$ is $G$-equivariant 
(cf.\ Corollary \ref{eqcon}) and it has an induced Hermitian metric $h^{0}$. The triple $(L^{0},\nabla^{0},h^{0})$ is 
a prequantization of the Hamiltonian action of $R_{G}$ on the smooth bundle of
symplectic manifolds $(I_{G}\backslash\mu^{-1}(0_{M})\to M,\tilde{\omega}^{0}$. 

Moreover, the K\"{a}hler structure on $S\to M$ induces a K\"{a}hler structure on $I_{G}\backslash\mu^{-1}(0_{M})\to M$.
\begin{theorem}\label{qure}
If $G$ is a proper regular Lie groupoid, then quantization commutes with internal reduction, i.e.\ there exists an 
isomorphism of vector bundles
\[Q(I_{G}\backslash\mu^{-1}(0_{M}),\omega^{0})\overset{\simeq}{\to}(Q(S,\tilde{\omega}))^{I_{G}}.\]
\end{theorem}
\begin{remark}
Note that for each $m\in M$ one can restrict the action of $G$ on $J:S\to M$ to a 
Hamiltonian action of the isotropy Lie group $G_{m}^{m}$ on $(S_{m},\omega|_{S_{m}})$. Likewise, the momentum map, 
prequantization data, K\"{a}hler structure all restrict to $S_{m}$, and hence give rise to a quantization commutes with
reduction statement as in the theorem for $G_{m}^{m}$, which is compact since $G$ is proper. 
This theorem was first formulated and proven for compact Lie group action by Guillemin and Sternberg (cf.\ \cite{GS}) and also
goes under the name of ``Guillemin-Sternberg conjecture''.
It is proven, in a more general form using $\spin^{c}$-Dirac operators, for compact Lie groups by
Meinrenken (cf.\ \cite{Me}), Meinrenken and Sjamaar (cf.\ \cite{MS}), Tian and Zhang (cf.\ \cite{TZ}) and Paradan (cf.\ \cite{Par}). For certain non-compact groups it 
is proved, in a somewhat different form using K-theory and K-homology (cf.\ Remark \ref{KKstuff}), 
by Hochs and Landsman (cf.\ \cite{HL}). 
For families the theorem was proven in \cite{Zha}, also within the setting of $\spin^{c}$-Dirac operators and K-theory.
\end{remark}
\begin{remark}
If the theorem holds, then the following diagram ``commutes'':
\[\xymatrix{Q(S,\tilde{\omega})\ar@{|->}[r]^-{R}&Q(I_{G}\backslash\mu^{-1}(0_{M})),\omega^{0})\simeq (Q(S,\tilde{\omega}))^{I_{G}}\\
(S,\tilde{\omega})\ar@{|->}[r]_-{R}\ar@{|->}[u]^{Q}&((I_{G}\backslash\mu^{-1}(0_{M})),\omega^{0})\ar@{|->}[u]_{Q},}\]
where $R$ denotes symplectic and quantum reduction and $Q$ denotes quantization. One sometimes abbreviates the theorem
by writing $[Q,R]=0$.
\end{remark}
\begin{proof}
We shall construct a morphism
\[(Q(S,\tilde{\omega}))^{I_{G}}\to Q(I_{G}\backslash\mu^{-1}(0_{M}),\omega^{0}).\]
The inclusion
$i:\mu^{-1}(0_{M})\hookrightarrow S$ induces a map
\[i^{*}:\Gamma^{\infty}(L)^{I_{G}}\to \Gamma^{\infty}(L|_{\mu^{-1}(0_{M})})^{I_{G}},\]
where the superscript $I^{G}$ means that we restrict to equivariant sections. These are the sections fixed under the action of
$I_{G}$ on $\Gamma^{\infty}(L)$ as a bundle over $M$.
Moreover, the quotient map $p:\mu^{-1}(0_{M})\to \mu^{-1}(0_{M})/I_{G}$ induces an isomorphism
\[p^{*}:\Gamma^{\infty}(L^{0})\to \Gamma^{\infty}(L|_{\mu^{-1}(0_{M})})^{I_{G}}.\]
Because of the equivariance of the K\"{a}hler structure and the connection the composition $(p^{*})^{-1}\circ i^{*}$
induces a map
\[Q(S,\tilde{\omega}))^{I_{G}}\to Q(I_{G}\backslash\mu^{-1}(0_{M}),\omega^{0})),\]
which is the one we wanted to construct.
This map is an isomorphism on each fiber, since the isotropy groups of a proper groupoid are compact groups for which the theorem is well established. 
Hence $\Psi$ is a continuous isomorphism of vector bundles.
\end{proof}
\begin{example}\label{princ9}
Suppose $H$ is a Lie group acting in Hamiltonian fashion on a K\"{a}hler manifold $(S,h,j,\omega^{S})$. Suppose 
$\pi:P\to M$ is a principal $H$-bundle and $\tau\in\Gamma^{\infty}(\Wedge^{1}(P)\otimes\mathfrak{h}$ a connection
1-form. As discussed in previous examples there exist a $J$-presymplectic form $\tilde{\omega}$ on 
$P\times_{H}S\to M$ and the action of the gauge groupoid $P\times_{H}P\rra M$ on $(P\times_{H}S\to M,\tilde{\omega})$
is Hamiltonian. Suppose $(L\to S,\nabla^{L},g)$ is prequantization of the action of $H$ on $S$.
Let $Q_{S}$ denote the associated quantization. In Example \ref{princ8} we have seen that the prequantization and
quantization of the $H$-action on $S$ give rise to a prequantization $(P\times_{H}L,\nabla,g')$ and a quantization 
$Q\to M$ of the action of the gauge groupoid $P\times_{H}P\rra M$ on $(P\times_{H}S\to M,\tilde{\omega})$.
In Example \ref{MWprinc} we saw that the internal Marsden-Weinstein quotient is isomorphic to the 
trivial bundle $(M\times\mu^{-1}(0)/H\to M,\omega^{0})$. The quantization of this bundle obviously equals the 
vector bundle $M\times Q(\mu^{-1}(0)/H, j^{0})\to M$. The statement of Theorem \ref{qure} follows from $[Q,R]=0$ for $H$
plus the following observation 
\begin{lemma}
The internal quantum reduction of $Q(P\times_{H}S,j')$ is isomorphic to 
\[M\times Q_{S}^{H}.\]
\end{lemma}
\begin{proof}
In Example \ref{princ8} we proved that $Q(P\times_{H}S,j')\simeq P\times_{H}Q(S,j)$. 
An element $[p,\xi]\in P\times_{H}Q(S,j)$ is fixed under all $[p',p]\in I_{P\times_{H}P}$ whenever $\xi$ is
fixed under all $h\in H$. Hence the statement follows.
\end{proof}
\end{example}

Suppose $\pi:G\to U(E)$ is a unitary representation of $G$ on a 
vector bundle $E\to M$ with hermitian structure.
\begin{definition}
The \textbf{quantum reduction} of $\pi:G\to U(E)$ is the quotient vector bundle $E^{G}:=G\backslash E^{I_{G}}\to G\backslash M$ (these spaces are in general non-Hausdorff).
\end{definition}

\begin{corollary}\label{qure2}
If all conditions for Theorem \ref{qure} are satisfied, then quantization
commutes with symplectic reduction, i.e. there exists an isomorphism of continuous vector bundles
\[Q(G\backslash\mu^{-1}(0_{M}),\tilde{\omega}^{0})\overset{\simeq}{\to}(Q(S,\tilde{\omega}))^{G}\]
\end{corollary}

\begin{example}
As a (very) basic example we consider the pair groupoid $M\times M\rra M$ for a manifold $M$, acting on $(id:M\to M,0)$.
A momentum map is any map $\mu\in\Omega^{1}(M)$ such that $d\mu=0$. Since we assume $0_{M}\subset\im(\mu)$, $\mu$ has to be zero.
Hence the Marsden-Weinstein quotient equals $(*,0)$.
Recall that quantization and prequantization line bundles coincide in this case.
The only quantization representation of $T M$ that integrates to a representation of $M\times M$ is 
the trivial representation on the trivial complex line bundle $M\times\CC\to M$. Obviously, the quantum 
reduction of such a bundle is $\CC\to *$. The Marsden-Weinstein quotient $(*,0)$ is indeed quantized
by $\CC$.
\end{example}

\begin{example}
The previous example is a special case of gauge groupoid considered in the previous Examples \ref{princ}, \ref{princ2},
\ref{princ3}, \ref{princ4}, \ref{princ5}, \ref{princ8}, \ref{princ9}, where $P\to M$ is a principal $H$-bundle and $H$ 
a Lie group. One easily sees that the full quantum reduction of 
\[Q(P\times_{H}S,j')\simeq P\times_{H}Q(S,j)\]
is isomorphic to $Q(S, j)^{H}$. Moreover, the Marsden-Weinstein quotient equals $Q(\mu^{-1}(0)/H, j^{0})$. Hence, in this
example, it is particularly clear how $[Q,R]=0$ for gauge groupoids reduces to $[Q,R]=0$ for Lie groups.
\end{example}

\subsection{The orbit method}\label{orbitmethod}
To investigate and illustrate a possible orbit method we treat some examples in this section. We have not yet arrived at a full formulation.

A representation of Lie groupoid $G$ or Lie algebroid $\AA$ is said to be irreducible if
is has no proper $G$-(or $\AA$-)invariant subvector bundle.
\begin{example}
Suppose $H$ is a Lie group and $P\to M$ an $H$-principal bundle. 
The isotropy groupoid of the gauge groupoid
$P\times_{H}P$ is isomorphic to $P\times_{H}H$, where the action of $H$ on $H$ is by conjugation.
Indeed, define a map
\[I_{P\times_{H}P}\to P\times_{H}H\]
by
\[[p,q]\mapsto [p, q p^{-1}],\]
where $q p^{-1}$ is the unique element in $H$ such that $(q p^{-1})p=q$.
This is well-defined, since $[h p, h q]\mapsto [h p, h q p^{-1}h^{-1}]=[p, q p^{-1}]$.
Hence the bundle of Lie algebras $\AA(I_{P\times_{H}P})$ is isomorphic to $P\times_{H}\mathfrak{h}$, where the action
of $H$ on $\mathfrak{h}$ is the adjoint action. Moreover, the dual bundle $\AA(I_{P\times_{H}P})^{*}$ is isomorphic 
to $P\times_{H}\mathfrak{h}^{*}$, where the action of $H$ on $\mathfrak{h}$ is the coadjoint action.
From a coadjoint orbit $\mathcal{O}\subset\mathfrak{h}^{*}$ one can construct a bundle 
$P\times_{H}\mathcal{O}\subset P\times_{H}\mathfrak{h}^{*}$, which is easily seen to correspond to a coadjoint orbit
in $\AA(I_{P\times_{H}P})^{*}$.
\begin{lemma}
All coadjoint orbits of the gauge groupoid $P\times_{H}P$ are isomorphic to $P\times_{H}\mathcal{O}$ for a coadjoint
orbit $\mathcal{O}\subset\mathfrak{h}^{*}$.
\end{lemma}
If we choose a connection on $P$, then we can extend the symplectic structure on $\mathcal{O}$ to a 
$J$-presymplectic form, such that the action of the gauge groupoid is Hamiltonian.
In Proposition \ref{princ5} we proved that the quantization of this action does not depend on the choice of connection.
Actually, the representation of $P\times_{H}P$ is obtained from quantization of
the action of $H$ on $\mathcal{O}$ by Morita equivalence. In general, since $H$ and $P\times_{H}P$ 
are Morita equivalent, there exists a bijection between irreducible unitary representations of $H$ and 
irreducible unitary representations of $P\times_{H}P$. From the above lemma we see
that there is also a bijection between the coadjoint orbits of $H$ and $P\times_{H}P$. From this we conclude that for 
as far as the orbit method (using geometric quantization) works for $H$, it works for $P\times_{H}P$ too. 
\end{example}
\begin{remark}
A natural guess would be that, in general, there exists a correspondence between (isomorphism classes of) irreducible Hermitian Lie algebroid representations and smooth families of coadjoint orbits. However, such a correspondence almost always fails.  
The integration of Lie algebroid representations to Lie groupoid representations plays an important r\^{o}le, as it does in the Lie algebra/group case. For example consider the pair groupoid $M\times M$, with Lie algebroid $T M$. 
It has only one (trivial) smooth bundle of coadjoint orbits, namely $\AA^{*}(I_{G})=M\times\{0\}$.
Quantization gives, in general, many irreducible unitary representations of $T M$, namely flat connections on 
complex Hermitian line bundle over $M$. But, the only one that integrates to a representation 
of $M\times M$ is the trivial line bundle.
\end{remark}
\begin{example}\label{cex}
The previous example gives rise to an example where the orbit method for Lie groupoids fails. 
Consider the fundamental groupoid $\pi_{1}(M)$ of $M$. Its Lie algebroid is $T M$, hence there is only one coadjoint orbit, namely the zero orbit. Quantization can in general give rise to many non-isomorphic irreducible unitary representations of $T M$. This time the representations all integrate to unitary representations of $\pi_{1}(M)$ given by parallel transport. Hence, in this case, there is no bijection between smooth families 
of coadjoint orbits and
irreducible unitary representations of $\pi_{1}(M)$. This is reflected by the fact that one can choose different presymplectic
forms on the zero coadjoint orbit. Indeed, one has, in general, just one coadjoint orbit, but $\Pic_{T M}(id)=\Pic(M)=\check{H}^{2}(M,\ZZ)\not=0$.
\end{example}
\begin{remark}
Instead, one might hope for a correspondence between symplectic leaves the Poisson structure on $\AA^{*}$ and isomorphism classes of unitary representations of the groupoid. The same Example \ref{cex} shows this will not work. We think that a good formulation of an orbit method for groupoids should incorporate the K-theory of $M$, but we have to leave this as an open problem. The next two examples show that the orbit method at least works to some extend, when $M$ is contractible and hence its K-theory is zero, i.e.\ if there are no non-trivial vector bundles. 
\end{remark}
\begin{example}
Now we shall consider a non-regular groupoid. A simple example is given by the action groupoid $G:=S^{1}\ltimes\RR^{2}$ of the action of the
circle $S^{1}$ on the plane $\RR^{2}$ by rotation around the origin. 
The dual of the bundle of Lie algebras associated to the isotropy groupoid is given by 
\[(\AA^{*}(I_{G}))_{(x,y)}\simeq\left\{\begin{array}{ll}\RR&\mbox{ if }(x,y)=(0,0)\\
0&\mbox{ if }(x,y)\not=(0,0)
\end{array}\right.\]
The only {\em smooth} bundle of coadjoint orbits is the trivial one $\{0\}\times\RR^{2}\simeq\RR^{2}$. 
This suggests that all smooth representations of $G$ are trivial at the origin. Close inspection shows that this is indeed true.
The only irreducible representation of $G$ up to isomorphism is obtained by geometric quantization. Indeed, a $J$-presymplectic 
form is necessarily $0$. Hence a smooth momentum map $\RR^{2}\to\AA^{*}(G)$ lifting the zero inclusion 
$\RR^{2}\to \AA^{*}(I_{G})$ is given by $(r,\alpha)\mapsto(f(r,\alpha),(r,\alpha))$, where $f$ is smooth and $f(r,\alpha)=0$ if $r=0$
(using polar coordinates $(r,\alpha)$).
The prequantum line bundle $L$ is necessarily the trivial one (since $\RR^{2}$ is contractible) and the prequantization(=quantization) representation 
\[\AA(G)\simeq \RR\times\RR^{2}\to(\mathfrak{u}(1)\times\RR^{2})\oplus T \RR^{2}\simeq \mathcal{D}(L)\] 
is given by
\[(X,(r,\alpha))\mapsto (-2\pi i f(r,\alpha), X\frac{d}{d\alpha}),\]
which integrates to a representation
\[S^{1}\ltimes\RR^{2}\to U(L)\simeq \RR^{2}\times U(1)\times\RR^{2}\]
\[(\beta,(r,\alpha))\mapsto((r,\alpha+\beta),e^{2\pi i(f(r,\alpha+\beta)-f(r,\alpha))},(r,\alpha)).\]
\end{example}
\begin{example}
Even for continuous families of Lie groups geometric quantization and the orbit method work, although one should proceed
with caution. For example, consider the 2-sphere $S^{2}\subset\RR^{3}$. It can be seen as a continuous 
family of Lie groups under the projection $S^{2}\to [-1,1]$ given by $(x,y,z)\mapsto x$.
The dual of the associated bundle of Lie algebras is given by
\[(\AA^{*}(S^{2})_{x}\simeq\left\{\begin{array}{ll}\RR&\mbox{ if }x\in]-1,1[\\
0&\mbox{ if }x=\pm 1\end{array}\right.\]
The image of any continuous section $\theta:[-1,1]\to\AA^{*}(S^{2})$ is a continuous family of coadjoint orbits (which are points).
A momentum map is given by inclusion $\mu:\theta([-1,1])\hookrightarrow \AA^{*}(S^{2})$.
A prequantum line bundle is again necessarily trivial $L=[-1,1]\times\CC$.
The prequantum representation is given by
\[(x,X)\mapsto 2\pi i\la\mu,X\ra.\]
The remarkable feature of this example is that one can allow a $\theta$, and hence $\mu$, which is {\em not}
continuous at $x=\pm 1$, namely a fixed $k\in\ZZ$ on $]0,1[$ and $0$ in $\pm 1$, and still find a continuous representation after integration:
\[(x,\alpha)\mapsto e^{2\pi i \mu \alpha}.\]
This is a particular instance of the fact there exist non-continuous vector fields which still induce
homeomorphisms. Realizing this fact, an orbit method should allow 
families of coadjoint orbits that are non-continuous at the points $x=\pm 1$.
\end{example}

\end{document}